\documentclass[11pt,a4paper,twoside]{amsart}
\usepackage[utf8]{inputenc}
\usepackage[english]{babel}
\usepackage[T1]{fontenc}
\usepackage{ucs}
\usepackage{amsmath}
\usepackage{amsfonts}
\usepackage{amssymb}
\usepackage{amsthm}
\usepackage[all]{xy}
\usepackage{enumerate}
\usepackage{enumitem}
\usepackage{upgreek}
\usepackage{tipa}
\usepackage{tikz}
\usepackage{tikz-cd}
\usepackage{pgfplots}
\usetikzlibrary{patterns}
\usetikzlibrary{matrix,arrows}
\definecolor{refblue}{RGB}{0, 0, 153}
\usepackage[colorlinks,linkcolor=refblue,citecolor=citegreen,urlcolor=red,bookmarks=false,hypertexnames=true]{hyperref}
\usepackage[nameinlink]{cleveref}
\usepackage{ifthen}

\usepackage{verbatim} 
\usepackage{mathtools}
\usepackage{url}
\usepackage{mathrsfs}
\usepackage{indentfirst}
\usepackage{fancyhdr}
\usepackage{colonequals}
\usepackage{enumitem}

\usepackage[textwidth=20mm, textsize=tiny]{todonotes}
\usepackage[skip=5pt, indent=10pt]{parskip}

\definecolor{citegreen}{RGB}{0, 115, 0}
\hypersetup{colorlinks=true,linkcolor=refblue,citecolor=citegreen}

\usepackage[top=3.5cm, bottom=3.5cm, left=2.4cm, right=2.4 cm]{geometry}

\usepackage{thmtools, thm-restate}

\usepackage{appendix}

\newtheorem*{thm*}{Theorem}
\newtheorem*{conj*}{Conjecture}

\newtheorem{theorem}{Theorem}[section]

\newtheorem{lemma}[theorem]{Lemma}

\newtheorem{proposition}[theorem]{Proposition}

\newtheorem{corollary}[theorem]{Corollary}
\newtheorem{maintheorem}{Theorem}

\theoremstyle{definition}
\newtheorem{example}[theorem]{Example}
\newtheorem{remark}[theorem]{Remark}

\newtheorem{definition}[theorem]{Definition}
\newtheorem{notation}[theorem]{Notation}

\newcommand{\bbQ}{\mathbb Q}

\newcommand{\bbZ}{\mathbb Z}

\newcommand{\Spec}{\operatorname{Spec}} 
\newcommand{\End}{\operatorname{End}} 
\newcommand{\Hom}{\operatorname{Hom}} 
\newcommand{\id}{\operatorname{id}} 
\newcommand{\Pic}{\operatorname{Pic}} 
\newcommand{\pr}{\operatorname{pr}} 

\newcommand{\CH}{\operatorname{CH}} 

\newcommand{\ch}{\operatorname{ch}}
\newcommand{\Cer}{\operatorname{Cer}}
\newcommand{\Q}{\mathbb{Q}}
\newcommand{\intt}{\textup{int}}

\newcommand{\Corr}{\operatorname{Corr}}

\numberwithin{equation}{section}

\makeatletter

\renewcommand\subsubsection{\@secnumfont}{\bfseries}%
\renewcommand\subsubsection{\@startsection{subsubsection}{3}
  \z@{.5\linespacing\@plus.7\linespacing}{-.5em}%
  {\normalfont\bfseries}}

\makeatletter
\def\blfootnote{\gdef\@thefnmark{}\@footnotetext}
\makeatother

\title{On Modified Diagonal Cycles and the Beauville Decomposition of The Ceresa Cycle}

\author[L.~Lagarde]{Lucas Lagarde}
 \address{Lucas Lagarde, Institut Galilée, Université Sorbonne Paris Nord, 93430 Villetaneuse, France}
 \email{lagarde@math.univ-paris13.fr}

\author[M.~Moakher]{Mohamed Moakher} 
 \address{Mohamed Moakher, Department of Mathematics, University of Pittsburgh, Pittsburgh, USA}
 \email{mom224@pitt.edu}

\author[M.~Porzio]{Morena Porzio} 
\address{Morena Porzio, Department of Mathematics, University of Toronto, Toronto, ON, Canada}
\email{m.porzio@utoronto.ca}

\author[J.~Rawson]{James Rawson} 
\address{James Rawson, School of Mathematics and Statistics, University of Glasgow, Glasgow, UK}
\email{james.rawson@glasgow.ac.uk}

\author[F.~Trejos Suárez]{Fernando Trejos Suárez} 
 \address{Fernando Trejos, Department of Mathematics, Princeton University, Princeton, NJ, USA}
 \email{trejos@princeton.edu}

\pgfplotsset{compat=1.18}

\begin{document}

\begin{abstract}
Let $C$ be a curve of genus $g \geq 2$, and let $J$ be its Jacobian.
The choice of a degree 1 divisor $e$ on $C$ gives an embedding of $C$ into $J$; we denote by $[C]^{e}\in \mathrm{CH}( J;\mathbb{Q} ) $ the class in the Chow group of $J$ defined by its image. It is known that the vanishing of the Ceresa cycle $\mathrm{Cer}(C,e):=[C]^e - [-1]_* [C]^e$ is equivalent to both the vanishing of the 1st Beauville component $[C]_{(1)}^e$ and the vanishing of the 3rd Gross--Kudla--Schoen modified diagonal cycle $\Gamma^3(C,e) \in \mathrm{CH}(C^3;\mathbb{Q})$. We extend this result to show that the vanishing of the $s$-th Beauville component $[C]^e_{(s)}$ for $s \geq 1$ is equivalent to the vanishing of the $(s+2)$-nd modified diagonal cycle $\Gamma^{s + 2}(C, e) \in\mathrm{CH}(C^{s+2};\mathbb{Q})$. Moreover, we establish ``successive vanishing'' results for these cycles. We apply our results to study the rational (non)-triviality of $[C]_{(s)}^e$ in the special case $s=2$. Finally in the $s=1$ case, we show an integral refinement to the original statement, relating the order of torsion of $\mathrm{Cer}(C,e) \in \mathrm{CH}(J;\mathbb{Z})$ to that of $\Gamma^3(C,e) \in \mathrm{CH}(C^3;\mathbb{Z})$. 
\end{abstract}

\maketitle
\begin{center}
    \rule{150mm}{0.2mm}
\end{center}

\setcounter{tocdepth}{1}
\tableofcontents

\section{Introduction}
Let $k$ be an algebraically closed field of arbitrary characteristic, and let $C$ be a smooth connected projective curve of genus $g\geq 2$ over $k$. For any degree 1 divisor $e$ with integer coefficients, there is an embedding $\iota_e$ of $C$ into its Jacobian $J$, given by $x \mapsto \mathcal{O}_C(x - e)$. The image $[C]^e = \iota_{e, *}(C)$ is a natural class in the Chow group of the Jacobian, and is connected to interesting geometric invariants such as the Ceresa cycle, defined as $\Cer(C, e):= [C]^e - [-1]_* [C]^e$. This cycle is homologically trivial, but for a very general curve of genus $g \geq 3$, {it has infinite order modulo algebraic equivalence} (see \cite{Ceresa}, \cite{Fk96}). 
Since its discovery, determining the precise locus of curves $C$ for which the Ceresa cycle is non-vanishing, finding vanishing criteria, and constructing families of curves on which the Ceresa cycle does vanish, have become active areas of research (see for instance \cite{beauville-nonhyp}, \cite{beauville-schoen}, \cite{GZ},\cite{QZ24} \cite{LS23}, \cite{LS24}).

The behaviour of these cycles is typically studied in the Chow ring with $\mathbb{Q}$-coefficients (to suppress torsion), which we will simply write as $\CH(X)$ for a variety $X$ over $k$. 

As we are working with rational coefficients, we have the Beauville decomposition $$\CH(J) = \bigoplus_{i,s} \CH_{i, (s)}(J)$$ where $\CH_{i, (s)}(J) = \{c \in \CH_i(J) : [n]_* c = n^{2i + s} c, \forall n \in \mathbb{Z}\}$, for $[n]: J\rightarrow J$ the multiplication-by-$n$ map. For any class $z \in \CH_i(J)$, we will denote its component in $\CH_{i, (s)}(J)$ by $z_{(s)}$.
Applying the Beauville decomposition to the Ceresa cycle, we have $\mathrm{Cer}(C,e) =2\cdot {\sum}_{s \geq 1\ \text{odd}}[C]_{(s)}^e$. It therefore becomes interesting to study the Beauville components $[C]^e_{(s)}$ of the curve class in order to study the vanishing of the Ceresa cycle itself.

When working modulo algebraic equivalence, the classes $[C]^e_{(s)}$ are independent of the choice of $e$, and have been extensively studied. For example, E. Colombo and B. van Geemen \cite{CvG93} showed that if $C$ admits a morphism of degree $d$ to $\mathbb{P}^1_k$, then $[C]_{(s)}^e$ is algebraically trivial for all $s\geq d-1$. Conversely, they conjectured that for a very general curve (over $\mathbb{C}$) of genus $g\ge 2s+1$, the cycle $[C]_{(s)}^e$ is algebraically non-trivial (this is spelled out in \cite[Conjecture 1.4]{voisin13}); by the work of Voisin this is known at least for $g > s(s+3)/2$ \cite[Corollary 1.7]{voisin13}. 

There is another family of classes attached to the pair $(C, e)$, known as the modified diagonal or Gross--Kudla--Schoen cycles \cite{GK92}, \cite{GS95}. For any non-empty subset $I \subset \{1, ..., n\}$ we first consider the following cycle in $\CH_1(C^n)$,
\begin{equation}\label{eq-def delta}
\Delta_I^n(C, e) = \pr_{I}^*( \Delta^{(|I|)}_C) \cdot \prod_{i \notin I} \pr_{i}^* e ,
\end{equation} 
where $\pr_{I}: C^n \longrightarrow C^{|I|}$ is the projection onto components indexed by $I$, and $\Delta^{(m)}_C$ denotes the small diagonal in $C^m$. Then, the $n$-th modified diagonal cycle attached to $(C, e)$ is defined as
\begin{equation}\label{eq-gammandefintro}
\Gamma^n(C, e) = \sum_{\varnothing \neq I \subset \{1, ..., n\}} (-1)^{n - |I|} \Delta_I^n(C, e) \in \CH_1(C^n).
\end{equation}
These cycles are notable due to their applications to number theory: when $C$ is a Shimura curve, Gross and Kudla conjectured that $\Gamma^3(C,e)$ arises in special-value formulae for triple-product $L$-functions and their derivatives, which has been shown in various contexts (see \cite{YZZ}, \cite{DR}).

The relation between the modified diagonal cycles and the Ceresa cycle was first noticed by Colombo and van Geemen. 
Let $\sigma_{n}: C^n \rightarrow J$ be the natural extension of the map $\iota_e$ by the summation map.
Then \cite[Proposition 2.9]{CvG93} states that the two cycles $\sigma_{3,*}(\Gamma^3(C,e))$ and $3\cdot \Cer(C,e) $ in $\CH(J)$ are Abel--Jacobi equivalent, i.e., they are homologically trivial and have the same image under all Abel--Jacobi maps.
S.-W. Zhang later strengthened this result, and in particular established an equivalence between the vanishing of the cycles $\Gamma^3(C,e)$, $\Cer(C,e)$, and $[C]_{(1)}^e$, working modulo rational equivalence.

\begin{theorem}[{\cite[Theorem 1.5.5]{Zhang10}}]\label{Zhangtheorem}
Let $e$ be any degree 1 divisor on the curve $C$. {Then there is an equivalence between vanishing of cycles as follows:} 
\[
\Cer(C, e) = 0\ \  \iff \ \  [C]^e_{(1)} = 0 \ \ \iff \ \  \Gamma^3(C, e) = 0.
\]
Moreover, this vanishing can occur only if the equality $(2g - 2)e = K_C$ holds in $ \CH_0(C)$.
\end{theorem}

In light of this theorem, the choice of divisor $e=\xi:=K_C/(2g-2)$ has become standard for problems involving the Ceresa cycle. We note that Zhang's original motivation for considering $\xi$ is that this choice minimizes the Beilinson--Bloch height of the modified diagonal cycle $\Gamma^3(C,e)$ \cite[Section 1.3]{Zhang10}.

It is important to highlight that for a generic curve $C$ of genus $g \geq 4$, the cycle $\xi$ is not equivalent to a point in the group $\CH(C)$. This can be seen as a consequence of the generic non-vanishing of the Faber--Pandharipande cycle $(2g-2)^2\cdot \left(\xi\times \xi- \Delta_{C,*}(\xi)\right)\in \CH(C^2)$, which always vanishes when $\xi$ is equivalent to a point $x_0 \in C(k)$ (this cycle is introduced in \cite{GG03}; for the generic non-vanishing see \cite{GG03} and \cite{Yin15}). 
{For genus $g=3$, the same statement is true at least if we assume $k$ has characteristic zero (see \Cref{g3}).}

However, much work on (higher) Beauville components and (higher) modified diagonal cycles has been restricted to the case where the embedding $\iota_e$ is taken with respect to a basepoint $e=x_0$.
In particular, B. Moonen and Q. Yin \cite{MY16} have studied analogues of Zhang's results to higher Beauville components under these restrictions: they prove that the vanishing of $[C]^{x_0}_{(s)}$ is equivalent to the vanishing of $\Gamma^{s + 2}(C, x_0)$. Our first main theorem extends their work to an arbitrary choice of degree 1 divisor, $e$.
\begin{maintheorem}\label{maintheoremA}
    For $s \geq 1$ and any choice of degree 1 divisor $e\in \CH_0(C)$, we have 
    $$[C]^e_{(s)}=0\ \  \iff \ \  \Gamma^{s + 2}(C, e)=0.$$
\end{maintheorem}
This theorem is shown by studying the properties of the natural extension $\sigma_n\colon C^n \rightarrow J$ of the map $\iota_e$. The backwards implication is given by generalizing the method of Zhang and Colombo--van Geemen, which relates the cycle $\sigma_{s+2,*}(\Gamma^{s+2}(C,e))$ to $[C]^e_{(s)}$. For the forwards implication in \Cref{maintheoremA}, it proves convenient to replace $\Gamma^n(C,e)$ with another closely-related cycle, denoted $B^n(C,e) \in \CH_1(C^n)$. Its definition is inspired by a motivic interpretation of the modified diagonal cycles in \cite{MY16} (see \Cref{def:gammaandbeta} and \Cref{rmk-motivic}): according to this interpretation, $B^n(C,e)$ should be more closely related to the Jacobian (being the projection of $\Gamma^n(C,e)$ to the $\mathfrak{h}_1(C)^{\otimes n}$ component inside $\mathfrak{h}(C^n)$), and indeed this cycle arises more naturally while computing $\sigma_{n}^*([C]^e_{})$ (see \Cref{ytogamma}, \Cref{pullback of fourier transform}). We show that the vanishing of the cycle $B^n(C,e)$ is equivalent to the vanishing of $\Gamma^n(C,e)$, and we show that the vanishing of $[C]^e_{(n-2)}$ implies the vanishing of $B^n(C,e)$. We note $B^n(C,e)=\Gamma^n(C,e)$ when $e$ is a point, but the two differ in general.

Zhang's theorem \cite[Theorem 1.5.5]{Zhang10} also shows that the vanishing of $[C]^e_{(1)}$ must force the vanishing of $[C]^e_{(s)}$ for all $s>0$.
Restricting to the case where $e = x_0$, a more general statement can be deduced from the work of A. Polishchuk \cite{Polish05}, which shows that $[C]_{(s)}^{x_0} = 0$ implies $[C]_{(t)}^{x_0} = 0$ for all $t \geq s$. This has been shown separately for modified diagonals by K. O'Grady \cite{OGrady}, also restricting to the case where $e=x_0$. 
To extend Polishchuk and O'Grady's results to general $e$ with our methods, we need either a stronger assumption, or get a weaker conclusion.

\begin{maintheorem}\label{maintheoremB}
    For any $n \geq 3$, we have the following implication: 
    $$\Gamma^{n}(C, e) = \Gamma^{n+1}(C, e) = 0 \ \ \Longrightarrow \ \  \Gamma^k(C, e) = 0,\quad \text{for all } k \geq n.$$ Equivalently, for any $s \geq 1$, if $[C]^e_{(s )} = [C]^e_{(s+1)} = 0$, then $[C]^e_{(t)} = 0$ for all $t \geq s$. 
    
    If $n\geq 3$ and $\Gamma^n(C, e) = 0$, then we have at least $\Gamma^k(C, e) = 0$ for all $k \geq 2n - 2$.
\end{maintheorem}

We make note of two interesting zero-cycles which naturally arise in the proofs of \Cref{maintheoremA} and \Cref{maintheoremB}. 
Firstly, we consider the cycle
\begin{equation}
    \updelta^e := \iota_{e, *}(e) - [0] \in \CH_0(J).
\end{equation} 
This cycle is trivial if $e=x_0$, a basepoint, 
but can be non-trivial for general $e$; it loosely measures the defect in how the behaviour of $[C]^e$ deviates from that of $[C]^{x_0}$. 
One can show that $\updelta^e_{(0)} = \updelta^e_{(1)}=0$, but for $s \geq 2$, the appearance of $s$-th Beauville components $\updelta^e_{(s)}$ in many formulas is an obstacle to generalizing Polishchuk and O'Grady's results about the implications of the vanishing of $[C]^e_{(s)}$ and $\Gamma^{s+2}(C,e)$.
The second cycle of interest is denoted by $\gamma^{s}_e(e) \in \CH_0(C^s) $, where $s \geq 1$ (see \Cref{def:gammaandbeta}). {When $s = 2$, it is a multiple of the Faber--Pandharipande cycle.}
Relating these two cycles is key to our proof of \Cref{maintheoremA}. A corollary to our results is that, for $s \geq 2$, we have $\updelta^{e}_{(s)}=0$ if and only if $\gamma^{s}_e(e)=0$ (see \Cref{equivalentvanishingforzerocycles}). This is analogous to the relationship between $[C]^e_{(s)}$ and $\Gamma^{s+2}(C,e)$ in \Cref{maintheoremA}. Analogously to \Cref{maintheoremB}, we also obtain a successive vanishing result for $\gamma^s_e(e)$ and for $\updelta^e_{(s)}$ (see \Cref{successivegammavanishing}).

\begin{remark}
If $k=\overline{\bbQ}$, then 
the Beilinson conjectures \cite{Bei85} imply that the Albanese map $\CH_0(J)_{\deg=0}\to J(k)$ is injective. 
This would imply the vanishing of $\updelta^e$, as it lies in the kernel of this map (see \Cref{yinI2}). 
Furthermore, the filtration on $\CH(J)$ associated to the Beauville decomposition is conjectured to correspond to the (conjectural) Beilinson--Bloch filtration (see \cite[Conjecture 2.7]{voisin13}).
Such a correspondence would imply that if $k = \overline{\mathbb{Q}}$ one has $[C]_{(s)}^e=0$ for all $s\geq 2$. 
\end{remark}

We conclude with a closer investigation of our results in the cases of $s = 1, 2$, using the divisor $e=\xi=K_C/(2g-2)$. 
For $s=2$, combining our results with those of Yin \cite{Yin15}, we show the following result on the rational (non)-triviality of $[C]^\xi_{(2)}$:
\begin{maintheorem}\label{maintheoremc}
    For the generic curve $C$ of genus $g\ge 4$, we have $[C]_{(2)}^{\xi}\neq 0$. On the other hand, for any curve $C$ of genus $g\le 3$, we have $[C]_{(2)}^{\xi}=0$.
\end{maintheorem}

Since the gonality of a curve of genus $g$ (over an algebraically closed field) is always at most $\lfloor (g+3)/2\rfloor$, we note that $[C]_{(2)}^\xi$ is algebraically trivial for $g \leq 4$ by the aforementioned result of Colombo and van Geemen. Furthermore, their conjecture states that generically, $[C]_{(2)}^\xi$ is algebraically non-trivial for $g \geq 5$.

\begin{remark}\label{remarkaboutalgclosure}
Although we work over an algebraically closed field $k$, all of our results regarding the vanishing of cycles in Chow groups with rational coefficients hold also over an arbitrary base field, by the standard argument as in \cite[Lemma 1A.3]{blochbook}. 
For example in \Cref{maintheoremc}, the vanishing of $[C]_{(2)}^{\xi}$ for any curve of genus 3 over an algebraically closed field, extends to the generic curve of genus 3. 
\end{remark}

In the case $s = 1$, we recover the full statement of Zhang's theorem (see \Cref{sec:zhangnewproof}). This is a consequence of \Cref{maintheoremA}, \Cref{maintheoremB}, and \Cref{vanishing of gamma n-12}, which establishes that $\Gamma^3(C,e)=0$ implies $(2g-2)e=K_C$ (we note that our argument fills a gap in Zhang's proof, and clarifies a small error in the original theorem statement; see \Cref{qiuremark}).
Furthermore, we give an integral refinement of this result by establishing a relationship between the torsion order of the Ceresa cycle and that of the modified diagonal cycle, working with Chow groups with integral coefficients. 
To accomplish this, we combine our results with the decomposition of the integral Chow group $\CH(J;\bbZ)$ constructed in \cite[Theorem 4]{MP10}, which is related to the Beauville decomposition.

\begin{maintheorem}\label{integrality}
    Let $\xi_{\mathrm{int}}\in \CH_0(C;\bbZ)$ be an integral representative of $\xi$, and let $d\in \mathbb{Z}$. 
    \begin{enumerate}[label= \normalfont(\arabic*)]
        \item If $\Cer(C,\xi_{\mathrm{int}})\in \CH_1(J;\bbZ)[d]$, then $\Gamma^3(C,\xi_{\mathrm{int}})\in \CH_1(C^3;\bbZ)[2\cdot d]$.
        \item If $\Gamma^3(C,\xi_{\mathrm{int}})\in \CH_1(C^3;\bbZ)[d]$, then $\Cer(C,\xi_{\mathrm{int}})\in \CH_1(J;\bbZ)[M_{g+1}\cdot d]$.
    \end{enumerate}
    Here, $M_{g+1} = {\prod}_{\text{prime }p \leq g+1} p^{\ell_p},
$ with $\ell_p=\lfloor \frac{g}{p-1}\rfloor$ if $p\ge 3$, and $\ell_2=g-1$. 
\end{maintheorem}

\subsection*{Outline}
The paper is structured as follows. In \Cref{SectionCyclesonJac}, we develop some basic intersection theory results, in particular verifying that some standard facts about Chow rings with $\mathbb{Q}$-coefficients still hold with $\mathbb{Z}$-coefficients.
We also generalize our setting to the case that $e$ is an arbitrary degree 1 divisor with $\mathbb{Q}$-coefficients, rather than with $\mathbb{Z}$-coefficients. 
In \Cref{moddiagcycles}, we recall the definition of the modified diagonal cycles $\Gamma^n(C, e)$ and $\gamma^n_e(e)$; inspired by \cite{MY16}, we also define a motivic variant, $B^n(C, e)$, which represents a component of $\Gamma^n(C, e)$ under a K\"unneth decomposition. In the same section, we show the equivalence of the vanishing of $\Gamma^n(C, e)$ and $B^n(C, e)$, and establish the successive vanishing formulae. In \Cref{JandproductofC}, we show the equivalence of the vanishing of $[C]^e_{(s)}$ and $B^{s + 2}(C, e)$, by studying the behaviour of these classes under Fourier transforms and natural addition maps $\sigma_n\colon C^n \to J$. In \Cref{sec:applications}, we investigate an integral analogue of Zhang's result, and explore the behaviour of $[C]_{(2)}^e$ for general curves of genus $g \geq 4$.

\subsection*{Notation and Conventions}
Unless stated otherwise, all equivalences of algebraic cycles are in Chow rings $\CH(-)$ with $\Q $-coefficients, i.e., modulo rational equivalence up to torsion. When working with integral coefficients, we denote the Chow rings by $\CH(-;\bbZ)$.

We will work with an arbitrary degree 1 divisor $e\in \CH_0(C)$ in Sections \ref{SectionCyclesonJac} to \ref{JandproductofC}, and specialize to the case $e=\xi = K_C/(2g-2)$ for some applications in \Cref{sec:applications}.

Given two integers $i<j$, we will often write $[i,j]$ for the interval $\{i,i+1,\dots,j\}$.
If $I\subseteq \{1,\dots,n\}$ is a subset, let $\pr_{I} \colon C^n\longrightarrow C^{|I|}$ be the projection onto the factors indexed by $I$, and let $\widehat{\pr}_{I} \colon C^n\longrightarrow C^{n-|I|}$ for the projection onto the factors indexed by elements not contained in $I$.
For simplicity, we write $\pr_{i}\colon C^n\longrightarrow C$ for the case $I=\{i\}$ and $\pr_{ij}\colon C^n\longrightarrow C^2$ for the case $I=\{i,j\}$. Similarly, for any non-empty subset $I\subseteq \{1,\dots,n\}$, we will denote by $m_I\colon J^n\longrightarrow J$ the sum of the components indexed by $I$, and we write $m_{ij}$ for the case $I=\{i,j\}$.

For $n\ge 1$, we denote by $\Delta_C^{(n)}\colon C\to C^n$ the diagonal morphism $x\mapsto(x,\dots,x)$. 
For $n=2$, we will write $\Delta_C$ instead of $\Delta^{(2)}_C$. 
We will also denote by $\Delta^{(n)}_C$ the small diagonal, namely the image of the diagonal morphism.

\subsection*{Acknowledgements} 
We are very grateful to Ben Moonen for suggesting this project at Arizona Winter School 2024, for mentoring us during the workshop, and for providing us with valuable insights. 
We also wish to thank Jef Laga and Ziquan Yang for helpful conversations; Ari Shnidman for raising the question that led to \Cref{integrality}; and Congling Qiu for sharing the argument in \Cref{qiuremark}. We thank the anonymous referees for their careful reading of this article and for their helpful suggestions. 
Finally, we extend our thanks to the organizers of the Arizona Winter School, without whom we would not have had the opportunity to meet and work on this project.

JR was supported by the Warwick Mathematics Institute Centre for Doctoral Training during the preparation of this article, and gratefully acknowledges funding from the UK Engineering and Physical Sciences Research Council (Grant number: EP/W523793/1).

\section{Cycles on Jacobians}\label{SectionCyclesonJac}
\subsection{Embeddings of Curves}
Throughout, {we denote by} $C$ a smooth projective connected curve of genus $g\geq 2$ over an algebraically closed field $k$. 
Denote by $J= \Pic^0_C$ the Jacobian of the curve $C$, and by $J^t$ its dual abelian variety. We let $e$ be a degree $1$ divisor on $C$ with rational coefficients. Given a positive integer $n$, we define the map
\[
    a_n:C^n\longrightarrow \Pic^n_C,\quad (x_1,\dots , x_n)\mapsto \mathcal{O}_C(x_1+\dots + x_n).
\]
The map $a_1$ extends to an isomorphism $\CH_0(C;\bbZ)_{\deg=0}\xrightarrow{\sim} J(k)$. 
{Given $n,d\in \mathbb{Z}$, the tensor product of line bundles on $C$, $m\colon \Pic^n_C \times \Pic^d_C \longrightarrow \Pic^{n+d}_C$,
induces a morphism 
\[
 \CH(\Pic^n_C) \times \CH(\Pic^d_C)\longrightarrow \CH(\Pic^{n+d}_C),\quad  (x,y) \mapsto x\star y \colonequals m_*(x\times y)
\]
called the \emph{Pontryagin product}.
In particular for $d=n=0$, we recover the usual Pontryagin product on $\CH(J)$.} 

\begin{definition}\label{def:translation}
Given an integral divisor $\gamma$ of degree $d$ on $C$, the \emph{translation by} $\gamma$ is defined as the isomorphism 
\[
    t_{\gamma}: \Pic^n_C\longrightarrow \Pic^{n+d}_C, \quad \mathcal{O}_C(z) \mapsto \mathcal{O}_C(z+\gamma).
\]
At the level of Chow groups, the pushforward map $t_{\gamma,*}\colon \CH(\Pic^n_C;\bbZ)\to \CH(\Pic^{n+d}_C;\bbZ)$ coincides with the Pontryagin product ${[\gamma ]} \star (-)$.
\end{definition}

Since $J(k)$ is divisible (because $k=\overline{k}$), one can show that for any rational divisor $\gamma\in \CH_0(C)$ whose degree is an integer, there exists an integral divisor $\gamma_{\intt}\in \CH_0(C;\bbZ)$ representing it. Moreover, the choice of an integral representative is unique up to torsion. We now show that this choice does not affect the map $t_{\gamma,*}$ in the Chow group with $\bbQ$-coefficients.

\begin{proposition}\label{prop:independencepushforward}
Let $\gamma$ and $\gamma'$ be two integral divisors on $C$ of degree $d$. If $\tau=\gamma-\gamma'$ is torsion as an element of $\CH_0(C;\bbZ)$, then for every $z\in \CH(\Pic^{n}_C)$,
\[
    t_{\gamma, *}(z) = t_{\gamma', *}(z)\  \textup{ in} \ \CH(\Pic^{n+d}_C).
\]
\end{proposition}

\begin{proof}
Since $\tau$ is torsion, $n \tau = 0$ for some $n\in \mathbb{N}$. In particular, $[n]_* ({[\tau]} - [0])$ is zero in $\CH_0(J)$. Given that $[n]_*$ is injective on $\CH(J)$, we deduce that ${[\tau]} = [0]$. We then get the result from the equality $t_{\gamma,*}(z)=t_{\gamma',*}(z)\star {[\tau]}$.
\end{proof}

\begin{definition}\label{independenceofpushandpull}
{For a rational zero-cycle $\gamma$ of degree $d$, we define $t_{\gamma, *} : \CH(\Pic^n_C) \rightarrow \CH(\Pic^{n + d}_C)$ to be $t_{\gamma_\intt, *}$, for any integral representative $\gamma_\intt$. This is well-defined by the previous proposition.}
In particular, the pushforward map 
$$
\iota_{e, *}\colonequals t_{-e,*}\circ (a_1)_*\colon \CH(C) \longrightarrow \CH(J)
$$ 
is well-defined. 
Similarly, we define the pullback map $\iota_e^*$ as
\[
    \iota_e^*:\CH(J)\longrightarrow \CH(C),\quad \quad z\mapsto \pr_{C, *}(\pr_J^*z \cdot \Gamma_{e_{\mathrm{int}}}),
\]
where ${e_{\mathrm{int}}}$ is any integral representative of $e$, and $\Gamma_{e_{\mathrm{int}}} \in \CH(C \times J; \mathbb{Z})$ is the graph of $\iota_{e_{\mathrm{int}}}$. Indeed, by the same principle as before, different choices of ${e_{\mathrm{int}}}$ produce graphs that differ only by torsion. 
\end{definition}

\subsection{Recollection of the Beauville decomposition and Fourier transform}\label{subsectiononBeauvilleandFourier}
For further details on the results in this subsection, one may consult \cite{MoonenAWS}.

We consider the canonical principal polarization $\phi\colon J\xrightarrow{\sim} J^t$ of $J$. 
{Let $\mathscr{P}$ be the Poincaré bundle on $J\times J^t$, and  
let $\Theta$ be the theta divisor in $\Pic^{g-1}_C$.
The Mumford bundle is the line bundle on $J\times J$ defined as
\[
\Lambda_\Theta \coloneqq (\id\times \phi)^*\mathscr{P}= (\phi\times \id)^*\mathscr{P}^t,
\]
where $\mathscr{P}^t$ denotes the pullback of the Poincaré bundle along the swap-map $\sigma\colon J^t\times J\rightarrow J\times J^t$. 
It can be explicitly described as
\begin{equation}\label{eq:LambdaTheta}
    \Lambda_\Theta = \mathcal{O}_{J\times J}(m^*\Theta_{\gamma} - \pr_1^*\Theta_{\gamma}-\pr_2^*\Theta_{\gamma}),
\end{equation}
where $\Theta_\gamma$ is the translation of $\Theta$ by a divisor $\gamma$ of degree $g-1$. In particular, the right hand side of \eqref{eq:LambdaTheta} is independent of the choice of $\gamma$.}
We will frequently use the fact that 
for all $n\in \mathbb{Z}$, $$([n]\times \id)^*\Lambda_\Theta=(\id\times [n])^*\Lambda_\Theta=\Lambda_\Theta^{\otimes n},$$
which follows from the analogous statement on $\mathscr{P}$ (see \cite[Cor. 4.10]{MoonenAWS}).

Now, recall that the Chow ring $\CH(J)$ can be decomposed as a direct sum 
\[
    \CH(J) = \bigoplus_{i,s} \CH_{i,(s)}(J) = \bigoplus_{i,s} \CH^i_{(s)}(J),
\]
called the Beauville decomposition. Here, $\CH_{i, (s)}(J) = \{c \in \CH_i(J) \ | \ [n]_* c = n^{2i + s} c, \forall n \in \mathbb{Z}\}$, and 
the groups $\CH_{i,(s)}(J)=\CH^{g-i}_{(s)}(J)$ vanish outside the range $-i\le s \le g-i$
(see \cite{Beau83}).
Given $z\in \CH_i(J)$, denote by $z_{(s)}$ its $s$-th Beauville component in $\CH_{i, (s)}(J)$.

Beauville conjectured that $\CH_{i,(s)}(J)=0$ whenever $s < 0$, and proved this vanishing for $i\in \{0,1\}$ \cite[Proposition 3(a)]{Beau83}. In particular, we can write
$$[C]^e= [C]^e_{(0)} + ... + [C]^e_{(g - 1)},$$
and for any $z\in \CH_0(J)$, we have  $z=z_{(0)}+\cdots +z_{(g)}$.

We recall the existence of two maps from the group $\CH_{0}(J)$: the \emph{degree map} {$\deg \colon \CH_0(J)\to \bbQ$}, and the \emph{Albanese map}  
$$
\mathrm{Alb}\colon \CH_0(J)\to J(k)\otimes \bbQ, \quad \sum n_i [x_i] \mapsto \sum n_i x_i,
$$
where the second summation uses the group law on $J$. Let $I$ denote the kernel of $\deg$. Then $I = \oplus_{s \geq 1} \CH_{0, (s)}(J)$, and it forms an ideal with respect to the Pontryagin product. 
The kernel of $\mathrm{Alb}$ is then given by $I^{\star 2} = \oplus_{s \geq 2} \CH_{0, (s)}(J)$. Note that in general, we have that $I^{\star n} = \oplus_{s \geq n} \CH_{0, (s)}(J)$. 

The Chow group also admits an invertible Fourier transform and its dual 
$$\mathscr{F} \colonequals \ch(\mathscr{P})_* : \CH(J) \xrightarrow{\sim} \CH(J^t),\quad z\mapsto \pr_{J^t,*}\left( \ch(\mathscr{P})\cdot \pr_J^*(z) \right),
$$
$$\mathscr{F}^t \colonequals \ch(\mathscr{P}^t)_* : \CH(J^t) \xrightarrow{\sim} \CH(J), \quad z\mapsto \pr_{J,*}\left( \ch(\mathscr{P}^t)\cdot \pr_{J^t}^*(z) \right).
$$
Here, $\ch$ denotes the Chern character. These transforms satisfy $\mathscr{F}^t\circ \mathscr{F}=(-1)^g\cdot [-1]_*$. Moreover, $\mathscr{F}$ and $\mathscr{F}^t$ restrict to isomorphisms $\CH_{i,(s)}(J)\xrightarrow{\sim} \CH_{g-i-s,(s)}(J^t)$.

\subsection{Preliminary intersection theory results}\label{subsectinoonthetadiv}
\begin{notation}
    We denote by $\kappa$ a theta characteristic of $C$, that is, a cycle in $\CH_0(C;\bbZ)$ satisfying $2\kappa=K_C $. 
    Its associated Theta divisor $\Theta_{\kappa}$ is symmetric. We also fix an integral representative $\varepsilon=e_{\mathrm{int}}$ for $e$.
\end{notation}
We state the following two lemmas with integral coefficients rather than rational ones, as they will be needed for our integral results in \Cref{sec:integral}. 
\begin{lemma}\label{lemmapulbackthetatoC}
    Given $\gamma\in \CH_0(C;\bbZ)$ of degree $g-1$, we have the following equality in $\CH_0(C;\bbZ)$:
    \[
         \iota_{{\varepsilon}}^*(\Theta_{\gamma}) = K_C - \gamma + {\varepsilon}.
    \]
    In particular, $\iota_{\varepsilon}^*(\Theta_{\kappa})= \kappa + \varepsilon$.
\end{lemma}

\begin{proof}
    By \cite[Cor. 14.21]{MEvdGbook}, for a degree $g-2$ line bundle $\mathscr{L}$ on $C$, we have 
    \[
        (t_{\mathscr{L}}\circ a_1)^*\mathcal{O}_{\Pic^{g-1}_C}(\Theta) \cong \mathcal{O}_C(K_C)\otimes \mathscr{L}^{-1}.
    \]
    Setting $\mathscr{L}=\mathcal{O}_C(\gamma - \varepsilon)$, we get 
    \[(t_{\gamma-\varepsilon}\circ a_1)^* \Theta=
    (t_{-\varepsilon}\circ a_1)^*(t_{\gamma}^*\Theta)  = \iota_{\varepsilon}^*(\Theta_{\gamma}) =K_C - \gamma + \varepsilon. \qedhere
    \]
\end{proof}
The following identity is well known in the case that $e$ is a $k$-point. We verify it in the general case.
\begin{lemma}\label{lem-pullbackMB}
The following equality holds in $\CH(C^2;\bbZ)$: 
\begin{equation*}
   -(\iota_{\varepsilon}\times \iota_{\varepsilon})^*c_1(\Lambda_\Theta)
   =
   \Delta_C-(\varepsilon\times C)-(C\times \varepsilon),
\end{equation*}
where $c_1(\Lambda_\Theta)$ is the first Chern class of $\Lambda_\Theta$.
\end{lemma}
We note that the cycle $\Delta_C-(\varepsilon\times C)-(C\times \varepsilon)$ is precisely the idempotent correspondence $\pi_1$ arising in \Cref{sec-defofmod}.
\begin{proof}
Let $x_0\in C(k)$, and let $(\mathscr{P}_C,\alpha_C)$ be the universal line bundle of degree $0$ on $C\times J$ rigidified at $\{x_0\}\times J$. We first show that
\begin{equation}\label{eqc1Lambda|CtimesJ}
    -(\iota_{\varepsilon}\times \id_J)^* c_1(\Lambda_\Theta ) 
    =  
    c_1(\mathscr{P}_C) + \pr_J^* ( -\Theta_{(g-1)\varepsilon} + \Theta_{g\varepsilon-x_0}).
\end{equation}
To do this, consider the map
\[
f:C\times \Pic^g_C\longrightarrow \Pic^{g-1}_C,\quad \quad (x, \mathscr{L})\mapsto \mathscr{L}(-x).
\]
By \cite[Prop. 14.20]{MEvdGbook}, the line bundle $f^*(\mathcal{O}_{\Pic^{g-1}_C}(\Theta))\otimes \pr^*_{\Pic^g_C} \mathcal{O}_{\Pic^g_C}(-t^*_{-x_0}\Theta)$, with a fixed trivialization (the identity) at $\{x_0\}\times \Pic^{g}_C$, is the universal line bundle of degree $g$ on $C\times \Pic^g_C$.
So by the See-Saw Principle, the pullback
\[
(\id_C\times t_{g\cdot \varepsilon})^* [ f^*(\mathcal{O}_{\Pic^{g-1}_C}(\Theta))\otimes \pr^*_{\Pic^g_C} \mathcal{O}_{\Pic^g_C}(-t^*_{-x_0}\Theta) ]
\]
is isomorphic to $\mathscr{P}_C\otimes \pr_C^* \mathcal{O}_C(g\cdot \varepsilon)$. Since $f$ fits in the commutative diagram
\begin{center}
\begin{tikzcd}
    C\times J 
    \arrow[rr,"{m \circ ([-1]\circ\iota_{\varepsilon}\times \id_J)}"] 
    \arrow[d, swap, "{\id_C\times t_{g\cdot \varepsilon}}"] & 
    & 
    J 
    \arrow[d, "{t_{(g-1)\varepsilon}}"]
    \\ 
    C\times \Pic_C^g \arrow[rr,"f"] 
    &
    &
    \Pic_C^{g-1},
\end{tikzcd}
\end{center}
we can write
\begin{align*}
c_1(\mathscr{P}_C)+\pr_C^*(g\cdot\varepsilon)
&= \bigl(m\circ ([-1]\circ\iota_\varepsilon\times\id_J)\bigr)^* t_{(g-1)\varepsilon}^*\Theta
- (\id_C\times t_{g \varepsilon})^*\pr_{\Pic_C^g}^*(t_{-x_0}^*\Theta)\\
&= ([-1]\circ\iota_{\varepsilon}\times\id_J)^* m^*\Theta_{(g-1)\varepsilon}
- \pr_J^*(\Theta_{g \varepsilon-x_0}).
\end{align*}
By definition of $\Lambda_{\Theta}$, we have that
\[
    m^*\Theta_{(g-1)\varepsilon} = \pr_1^*\Theta_{(g-1)\varepsilon} + \pr_2^* \Theta_{(g-1)\varepsilon} + c_1(\Lambda_{\Theta}).
\]
So the term $([-1]\circ \iota_{\varepsilon}\times \id_J)^* m^*\Theta_{(g-1)\varepsilon}$ can be rewritten as
\[
\iota_{\varepsilon}^* [-1]^* \Theta_{(g-1)\varepsilon} + \pr_J^*(\Theta_{(g-1)\varepsilon}) + (\iota_{\varepsilon}\times \id_J)^* ([-1]\times \id_J)^* c_1(\Lambda_{\Theta}).
\]
Using the equalities $([-1] \times \id_J)^* \Lambda_\Theta \cong \Lambda_\Theta^{-1}$ and $[-1]^* \Theta_{(g-1)\varepsilon} = \Theta_{K_C -(g-1)\varepsilon}$, we find that
\[
- (\iota_{\varepsilon}\times \id_J)^* c_1(\Lambda_{\Theta}) = 
c_1(\mathscr{P}_C) + \pr_J^*(-\Theta_{(g-1)\varepsilon} + \Theta_{g\varepsilon-x_0} )
+
\pr_C^*(g\cdot \varepsilon - \iota_{\varepsilon}^* \Theta_{K_C -(g-1)\varepsilon}).
\]
By \Cref{lemmapulbackthetatoC}, we have that $g\cdot \varepsilon - \iota_{\varepsilon}^* \Theta_{K_C -(g-1)\varepsilon} =0$. Hence equation \eqref{eqc1Lambda|CtimesJ} follows.

Now in order to compute $-(\iota_{\varepsilon}\times \iota_{\varepsilon})^*c_1(\Lambda_\Theta)$, we use the See-Saw Principle {(see e.g. \cite[Section 2.1]{MoonenAWS})}. 
By universality of $\mathscr{P}_C$, the line bundle  $\mathscr{P}_C|_{C\times \iota_{\varepsilon}(x)}$ corresponds to the $k$-point $\mathcal{O}_{C}(x - \varepsilon)$ for any $x\in C(k)$.
Then by equation \eqref{eqc1Lambda|CtimesJ} and \Cref{lemmapulbackthetatoC}, the restriction to $C\times \{x\}$ of $-(\iota_{\varepsilon}\times \iota_{\varepsilon})^*c_1(\Lambda_{\Theta})$ is 
\begin{equation*}
    -(\iota_{\varepsilon} \times \iota_{\varepsilon})^* c_1(\Lambda_{\Theta})|_{C\times{x}}= c_1(\mathscr{P}_C|_{C\times \iota_{\varepsilon}(x)}) + \pr_2^*(x_0-\varepsilon)|_{C\times \{x\}} = x-\varepsilon.
\end{equation*}
On the other hand
\begin{equation*}
    (\Delta_C- \varepsilon \times C-C\times \varepsilon)|_{C\times x}\cong 
    x - \varepsilon\times x - 0 \cong x-\varepsilon.
\end{equation*}
Moreover, the restriction of $-(\iota_{\varepsilon}\times \iota_{\varepsilon})^*c_1(\Lambda_{\Theta})$ to $\{x_0\}\times C$ is
\begin{align*}
    c_1(\mathscr{P}_C|_{\{x_0\}\times C} ) + \pr_2^*(x_0-\varepsilon)|_{\{x_0\} \times C} = x_0 -\varepsilon,
\end{align*}
and the restriction of $(\Delta_C- \varepsilon \times C-C\times \varepsilon)|_{\{x_0\} \times C}$ is $x_0 -\varepsilon$ as well.
We conclude by another application of the See-Saw principle.
\end{proof}
\begin{definition}
    We define 
    $$\updelta^e:=\iota_{e, *}(e)-[0]\in \CH_0(J).$$
    In particular, if  $e=x_0$ is a $k$-point, then $\updelta^{x_0}=0$.
\end{definition}

\begin{remark}\label{yinI2}
The cycle $\updelta^e$ is actually in the ideal $I^{\star 2}$.
Indeed, write $e=\sum_i m_ix_i$ with  $m_i\in \mathbb{Q}$ and $x_i\in C(k)$, then
$ \iota_{e, *}(e) = \sum_i m_i [\mathcal{O}_C(x_i - e)].
$ In particular, $\iota_{e, *}(e)$ has degree $1$, so $\iota_{e, *}(e) - [0] = \updelta^e\in I$. 
Moreover,
$$  \mathrm{Alb}(\updelta^e) = \mathrm{Alb}\big( \sum_i m_i [\mathcal{O}_C(x_i - e)] - [0] \big) = \mathcal{O}_C\big(\sum_i m_i (x_i -e) \big) = \mathcal{O}_C(e - e)= 0. $$
Hence $\updelta^e\in \ker(\mathrm{Alb})=I^{\star 2}$, and so in particular $\updelta^e_{(0)}=\updelta^e_{(1)}=0$.
\end{remark}

As discussed in the introduction, a natural choice of degree-one divisor for questions involving the Ceresa cycle is $\xi=K_C/(2g-2)$. For a generic curve of genus $g\ge 3$, the cycle $\xi$ is not equivalent to a point. The proof of this fact in the case $g\ge 4$ is straightforward due to the generic non-vanishing of the Faber--Pandharipande cycle, and is provided in the introduction. We record here the proof for $g=3$ in characteristic zero. 
\begin{lemma}\label{g3}
    Assume the field $k$ has characteristic zero. If $C$ is the generic curve of genus $3$, then $\xi$ is not equivalent to a point in $\CH_0(C)$.
\end{lemma}
\begin{proof}
Suppose by way of contradiction that $C$ admits a point $x_0$ such that $N(4x_0-K_C)$ is principal for some positive integer $N$. By working locally around the Fermat curve $F_4$, we may assume that there exists a valuation ring $V$ and a curve over $\Spec(V)$ whose generic fibre is $C$ and special fibre is $F_4$ (see \cite[\href{https://stacks.math.columbia.edu/tag/00IA}{Tag 00IA}]{stacks-project}).
We map $C$ to its Jacobian by $x \mapsto [4x - K_C]$, and specialize to $F_4$. By the valuative criterion for properness, $x_0$ specializes to a point $x_0'$ on $F_4$.
Recall that any cusp $c_0$ of $F_4$ is an incident point of a hyperflex, so it satisfies $4c_0 = K_C$. 
Because $4x_0' - K_C$ is torsion, so is $x_0' - c_0$. 
Hence, by \cite[Theorem A]{coleman}, $x_0'$ is a cusp as well and $4x_0'=K_C$.

Since we are working in characteristic zero, the $N$-torsion subgroup scheme of the universal Jacobian is finite étale over the moduli space of genus $3$ curves. 
We may assume that $V$ is Henselian, and since specialization of finite \'etale covers is an equivalence of categories over Henselian pairs (see \cite[\href{https://stacks.math.columbia.edu/tag/0BQC}{Tag 0BQC}]{stacks-project}), $N$-torsion is injective under specialization. Hence if the generic curve meets the torsion locus, it meets the origin, and so $C$ has a hyperflex. By \cite[Ch 1., Proposition 4.9]{Verm83}, the hyperflex locus is a divisor of the moduli
space of smooth genus three curves, a contradiction.
\end{proof}

\section{Modified Diagonal Cycles}\label{moddiagcycles}
In this section, we recall the definition of the modified diagonal cycles $\Gamma^n(C,e)$ in the motivic formulation of Moonen--Yin \cite{MY16}. These can be viewed as a special case (for $z=[C]$) of the more general cycles $\gamma^n_e(z)$. This perspective naturally leads to the introduction of another family of cycles, denoted $\beta^n_e(z)$, arising as distinguished components of $\gamma^n_e(z)$ under a motivic Künneth decomposition. In fact $\beta^n_e(z)$ only differs from $\gamma^n_e(z)$ when $z=[C]$, in which case we denote it by $B^n(C,e)$. This cycle will play a central role, as we will later show in \Cref{JandproductofC} that it is directly related to the curve class $[C]_{(n-2)}^e$. 
 
We establish in \Cref{thm:gammavanishiffB} the equivalence between the vanishing of $\Gamma^n(C,e)$ and  of $B^n(C,e)$. We also prove successive vanishing results for the $\gamma^n_e$'s and the $\beta^n_e$'s in \Cref{successivegammavanishing} and \Cref{corollaryonvanishingofgamma}, which will serve as important ingredients in recovering Zhang's theorem together with its integral refinement in \Cref{sec:applications}.

Throughout this section, it will be useful to keep in mind the statement of the adjunction formula in this context. Namely, the equality
\begin{equation}\label{adjunction}
\Delta^*_C(\Delta_C)=-K_C
\end{equation} 
holds in $\CH_0(C)$. In particular this implies that $\Delta_C \cdot \Delta_C = -\Delta_{C,*}(K_C)$ in $ \CH_0(C^2)$.

\subsection{Definitions of Modified Diagonal Classes}\label{sec-defofmod}
We begin by recalling some aspects of the theory of motives that will be needed, mostly following \cite[Section 12]{MoonenAWS}. We then give the motivic description of the modified diagonal cycles as introduced by Moonen--Yin \cite{MY16}. 

For two smooth projective varieties $X$ and $Y$ over $k$, where $X$ is connected, we denote
\[
\mathrm{Corr}_i (X,Y) := \CH_{\dim(X)+i}(X \times Y)
\]
and $\mathrm{Corr}(X,Y) = \bigoplus_i \mathrm{Corr}_i (X,Y)$. The pushforward along such a correspondence $\rho\in \Corr(X,Y)$ is defined as
$$
    \rho_*\colon\CH(X) \longrightarrow 
    \CH(Y), \quad z  \mapsto \rho_*(z)\colonequals \pr_{Y, *}(\pr_X^*(z)\cdot \rho).
$$
If $Z$ is a third smooth projective variety over $k$, the composition of $\alpha \in \mathrm{Corr}(X,Y)$ and $\beta \in \mathrm{Corr}(Y,Z)$ is given by
\begin{equation*}
    (\alpha, \beta) = \pr_{XZ,*}(\pr^*_{XY}(\alpha) \cdot \pr_{YZ}^*(\beta)) \in \Corr(X,Z),
\end{equation*}
where the maps $\pr_{XZ}, \pr_{XY}, \pr_{YZ}$ are the obvious projections from the triple product $X \times Y \times Z$.
For $\rho \in \Corr(X,X)$ and $\tau \in \Corr(Y,Y)$, we may define the tensor product {as}
$$ 
\rho \otimes \tau = \pr_{1 3}^*(\rho)\cdot \pr_{24}^*(\tau)\in \Corr(X\times Y, X\times Y).
$$ 

We recall the definition of the category $\mathrm{Mot}_k$ of covariant Chow motives over $k$.
The objects in this category are triples $(X,p,m)$, where $X$ is a smooth projective {$k$-}variety, $p \in \mathrm{Corr}_0(X,X)$ is an idempotent correspondence (also called a ``projector''), and $m \in \mathbb{Z}$. The morphisms from $(X,p,m)$ to $(Y,q,n)$ are correspondences of the form 
$q \circ \Corr_{m-n}(X,Y)\circ p $. In particular, the identity morphism of $(X,p,m)$ is $p \circ \Delta_X \circ p$. There is a covariant functor $\mathfrak{h}: \mathrm{SmProj}_k \rightarrow \mathrm{Mot}_k$ sending $X$ to $ (X, \Delta_X, 0)$ and 
$f: X \rightarrow Y$ to the class of its graph $[\Gamma_f]$. Given a Chow motive $M$, we define its Chow group $\CH_i(M)\colonequals\Hom_{\mathrm{Mot}_k}\big(\textbf{1}(i),M\big)$ for $i\ge 0$, where $\textbf{1}(i)=(\mathrm{Spec}(k),\id,i)$.
In particular if $M=\mathfrak{h}(X)$ then $\CH_i(M)=\CH_i(X)$.

We do not make much use of the motivic structure. The most important fact is that projectors give rise to a direct sum decomposition of Chow groups.
In particular for the curve $C$, we use the $0$-cycle $e$ of degree $1$ to define \emph{orthogonal projectors}
\begin{equation}\label{eq-pi def}
    \pi_0\colonequals C\times e, 
    \quad 
    \pi_2\colonequals e\times C, 
    \quad 
    \pi_1\colonequals \Delta_C-\pi_0-\pi_2   
\end{equation}
in $\CH_1(C^2) = \Corr_0(C,C)$. 
These give a decomposition of $\mathfrak{h}(C)$ 
$$\mathfrak{h}(C) \cong \mathfrak{h}_0(C)\oplus \mathfrak{h}_1(C)\oplus \mathfrak{h}_2(C),$$
where $\mathfrak{h}_i(C) \colonequals (C,\pi_i, 0)$. In particular, $\CH(C)=\CH(\mathfrak{h}(C))=\oplus_i \pi_{i,*}\CH(\mathfrak{h}(C))$.
We also set $\pi_+\colonequals \pi_1+\pi_2=\Delta_C-\pi_0$, and 
$\mathfrak{h}_+(C) \colonequals (C,\pi_+,0)= \mathfrak{h}_1(C)\oplus \mathfrak{h}_2(C)$. 
Similarly, there is a motivic K\"unneth decomposition of $\mathfrak{h}(C^n)$ which induces a direct sum decomposition of $\CH(C^n)$: of particular interest for our purposes are the direct summands $\CH(\mathfrak{h}_{+\ldots+}(C^n)) = (\pi_+^{\otimes n})_* \CH(C^n)$ and $\CH(\mathfrak{h}_{1\ldots1}(C^n)) = (\pi_1^{\otimes n})_* \CH(C^n)$. We note that $\CH(\mathfrak{h}_{1\ldots1}(C^n))$ is a direct summand of $\CH(\mathfrak{h}_{+\ldots+}(C^n))$.

We briefly recall a convenient property of tensor products of motives that will be of use later.
If we let $\rho\in \CH(C^n)$, $\tau\in \CH(C^m)$ and $\pi \in \CH(C^2)$, then by definition of the tensor of correspondences we have the equality
\begin{equation}\label{decompositionofpushforward}
(\pi^{\otimes (n+m)})_*(\rho\times \tau) = (\pi^{\otimes n})_*(\rho) \times (\pi^{\otimes m})_*(\tau).   
\end{equation}

We can now give the motivic definition of the cycles we will be interested in.
\begin{definition}\label{def:gammaandbeta}
    Given $n\geq 1$, we define the homomorphisms $\beta_e^n,\gamma_e^n \colon \CH(C)\to \CH(C^n)$ 
    by
    $$\beta_e^n\colonequals (\pi_1^{\otimes n})_*\circ \Delta_{C,*}^{(n)},\quad \gamma_e^n\colonequals (\pi_+^{\otimes n})_*\circ \Delta_{C,*}^{(n)}.$$ In particular, we define the following \textit{modified diagonal} cycles in $\CH_1(C^n)$,
    $$ B^n(C,e)\colonequals \beta_e^n([C]), \quad \quad \Gamma^n(C,e)\colonequals \gamma_e^n([C]). $$
The fact that this matches our previous definition of $\Gamma^n(C,e)$ in \eqref{eq-gammandefintro} is shown in \Cref{explicitexpression}.
\end{definition}
\begin{remark}\label{rmk-motivic}
We briefly recall the motivic interpretation of \cite{MY16} to motivate this definition. 
A conjecture of Murre \cite{Murre} predicts the existence of a Chow--K\"unneth decomposition for smooth projective varieties over $k$. That is, for any such $X$ of dimension $d$, there exist (non-unique) idempotent correspondences $\pi_0,\ldots,\pi_{2d} \in \CH_d(X \times X)$ with the property that for any Weil cohomology theory  $H\colon \mathrm{SmProj}_k^{\text{op}}\rightarrow \mathrm{GrVect}_K$ with coefficients in a field $K$, which therefore factors through the category of \emph{contravariant} Chow motives
    \[
    \begin{tikzcd}
        \mathrm{SmProj}_k^{\mathrm{op}}
        \arrow[r, "{}^t\mathfrak{h}"]
        \arrow[rd, "H"']
        &
        \mathrm{Mot}_k^{\mathrm{op}}
        \arrow[d, "\overline{H}"]
        \\
        &
        \mathrm{GrVect}_K,
    \end{tikzcd}
    \]
    we have $H^i(X)=\overline{H}((X,\pi_{2d-i},0))$ (see e.g. \cite[Section 12.7]{MoonenAWS}). 
    The conjecture is known in particular for $d=1$, and the idempotents $\pi_0,\pi_1,\pi_2$ can be taken as in \eqref{eq-pi def}. The key observation of Moonen--Yin is that the modified diagonal cycles $\Gamma^n(C,e)$ of Gross--Schoen as defined in \cite{GS95} by the formula \eqref{eq-gammandefintro} can be rewritten as a projection under the idempotent $\pi_+$ as in \Cref{def:gammaandbeta}. 
    With this in mind, it is natural to instead consider $B^n(C,e) := (\pi_1^{\otimes n})_*(\Delta^{(n)}_C) \in \CH_1(C^n)$. 
    {Indeed it is the projection of $\Delta_C^{(n)}$ to $\mathfrak{h}_1(C)^{\otimes n}$, and
    the theory of the motive $\mathfrak{h}_1(C)$ is essentially that of the Jacobian}
    (see also \cite[\S 3.2]{Scholl94}). Moonen--Yin only consider the case in which $e=x_0$ is a $k$-point, where the two cycles $\Gamma^n(C,e)$ and $B^n(C,e)$ are equal. In the general case, the cycles differ, and in fact the idempotent $\pi_1$ and thus the cycle $B^n(C,e)$ arise more naturally in computations, as is clear for instance from the proof of \Cref{ytogamma}. One explicit reason for the appearance of $\pi_1$ is the well-known identity of \Cref{lem-pullbackMB}. 
\end{remark}

\begin{example}\label{examples}
For $n=1$, we have
\begin{enumerate}
    \item[(1)] $B^1(C,e)=0$ and $\Gamma^1(C,e)=[C]$.
    \item[(2)] $\beta_e^1(z)=\gamma_e^1(z)=z-\deg(z)e$, for any 0-cycle $z \in \CH_0(C)$.
\end{enumerate}
For $n=2$, we have
    \begin{enumerate}
        \item[(1)] 
        $B^2(C,e)=\Gamma^2(C,e)=\Delta_C-(C\times e)-(e\times C).$
        \item[(2)]  $\beta_e^2(e)=\gamma_e^2(e)=\Delta_{C,*}(e)-(e\times e)$. This cycle coincides, up to a scalar, with the Faber--Pandharipande cycle (see \cite{GG03}). 
        Since $C$ is connected, it is algebraically trivial.
    \end{enumerate}
For $n=3$, we have
\begin{enumerate}
\item $\Gamma^3(C,e) =
\Delta_C^{(3)} - (e\times\Delta_C) -(\pr_2 ^*(e) \cdot \pr_{13}^*(\Delta_C)) - (\Delta_C\times e)  +(C\times e \times e)  +(e\times C\times e) +(e\times e \times C),$
\item  $B^3(C,e)=\Delta_C^{(3)} - (e\times\Delta_C) -(\pr_2 ^*(e) \cdot \pr_{13}^*(\Delta_C)) - (\Delta_C\times e)-(C\times \Delta_{C,*}(e))-\pr_{13}^*(\Delta_{C,*}(e))-(\Delta_{C,*}(e)\times C) $.
\end{enumerate}
\end{example}
Unwinding the definition of $\gamma_e^n$ and $\beta_e^n$, we obtain an explicit description of these operators.
\begin{lemma}\label{explicitexpression}
We have the following expressions.
\begin{enumerate}
    \item[\textup{(1)}] If $z\in \CH_0(C)$ and $n\ge 1$, then
    $$  \gamma_e^n(z)=\beta_e^n(z)=(-1)^n \deg(z)\cdot \prod_{j=1}^n \pr_j^*(e) + \sum_{k=1}^n (-1)^{n-k}\cdot \big( \sum_{|I|= k} \pr^*_I(\Delta_{C,*}^{(k)}(z))\cdot \prod_{j\not\in I} \pr^*_j (e) \big).$$
    \item[\textup{(2)}] If $n\ge 1$, then
    $$  \Gamma^n(C,e) = \sum_{k=1}^n (-1)^{n-k}\cdot \big( \sum_{|I|= k} \pr^*_I(\Delta_{C}^{(k)})\cdot \prod_{j\not\in I} \pr^*_j (e) \big) = \sum_{\varnothing \neq I \subseteq \{1,\dots n \} }
    (-1)^{n-|I|}\cdot \Delta^n_I(C, e),$$
    where $\Delta^n_I(C,e)$ is as in \eqref{eq-def delta}.
    \item[\textup{(3)}] If $n\ge 3$, then
    $$ B^n(C,e) = \Gamma^n(C,e)-\sum_{i=1}^n \widehat{\pr}_i^*(\gamma_e^{n-1}(e)).$$
\end{enumerate}
\end{lemma}
\begin{proof}
    This is a direct computation.
    By definition, $\gamma_e^n(z)$ is equal to
    \begin{align*}
        &\pr_{[n+1, 2n], *}
        \big( \pr_{[1,n]}^*\Delta_{C,*}^{(n)}(z) \cdot \prod_{i=1}^n \pr_{i\, n+i}^*(\Delta_C - C\times e) \big) 
        \\
        &=\sum_{I\subseteq \{1,\dots,n\}} (-1)^{n-|I|}\cdot \pr_{[n+1, 2n], *} \big( \pr_{[1, n]}^* \Delta_{C,*}^{(n)}(z) \cdot \prod_{i\in I} \pr_{i\, n+i}^*(\Delta_C) \cdot \prod_{i\notin I} \pr_{n+i}^*(e) \big).
    \end{align*}
    Moreover, if $z$ is a $0$-cycle, then $\pr_{[1, n]}^*\Delta_{C, *}^{(n)}(z) \cdot \pr_{i\, n+i}^*(\pi_2)=0$, so the projection along $\pi_+^{\otimes n}$ coincides with the one along $\pi_1^{\otimes n}$. This yields the first formula. 
    
    If $z=[C]$, then the term corresponding to the case $I=\varnothing$ vanishes for dimension reasons, implying the second formula.
    Furthermore, by definition we have
    \begin{equation}\label{equationforGammaandB}
        B^n(C,e)
        = \Gamma^n(C,e) - \sum_{i=1}^n \pr_{[n+1, 2n], *} \big( \pr_{[1,n]}^*(\Delta_{C}^{(n)}) \cdot \pr_{i\, n+i}^*(\pi_2) \cdot \prod_{j \neq i} \pr_{n+j}^*(\pi_+) \big).
    \end{equation}
    Given the identity $\Delta_C^{(n)}\cdot \pr_i^*(e)=\Delta_{C, *}^{(n)}(e)$, the second term on the RHS in \eqref{equationforGammaandB} is equal to $\sum_{i=1}^n \widehat{\pr}_i^*(\gamma_e^{n-1}(e))$. 
    This implies the third formula.
\end{proof}


\subsection{Vanishing Results}\label{vanish_result}
We now study the consequences of the vanishing of $B^n(C, e)$.

Since $B^n(C,e)$ is the projection of $\Delta^{(n)}_C$ onto the component $\mathfrak{h}_{1\ldots1}(C^n)$, 
which is contained in $\mathfrak{h}_{+\ldots+}(C^n)$, we have $\pi_{1, *}^{\otimes n}\, ( \Gamma^n(C, e) ) = B^n(C, e)$. 
In particular, the vanishing of $\Gamma^n(C, e)$ implies the vanishing of $B^n(C, e)$.
In this subsection, we prove that the converse is true as well. 
In order to do so, we first establish two consequences of the hypothesis $B^n(C, e) = 0$ (\Cref{vanishing of gamma n-12} and \Cref{vanish3}), and then use these results to show that if $B^n(C, e)=0$ then $\Gamma^n(C, e) = 0$ as well (\Cref{thm:gammavanishiffB}).

{The essence of the following proofs is to use $B^n(C, e)$ as a correspondence between $C$ and $C^{n - 1}$, and as a correspondence between $C^2$ and $C^{n-2}$. 
These calculations could be done using the explicit formulae in \Cref{explicitexpression}, but this becomes cumbersome. Instead, we use the projectors to simplify the book-keeping at the cost of the argument being slightly less transparent.}

\begin{proposition}\label{vanishing of gamma n-12}
     Let $n\ge 3$, and assume that $B^n(C,e)=0$. Then, 
     $$\gamma^{n-2}_e(K_C+2e)=0, \quad \quad \text{and} \quad \quad \gamma^{n-1}_e(K_C+4e)=0.$$
\end{proposition}
\begin{proof}
     In this proof we will consider projections from $C^{2n}$, $C^n$ and $C^{2n-2}$. 
    To avoid confusion, we reserve $\pr_I:C^{2n}\longrightarrow C^{|I|}$ for the $I$-projection from $C^{2n}$, and denote by $\overline{\pr}_I:C^{n}\longrightarrow C^{|I|}$ the one from $C^{n}$ and by $\widetilde{\pr}_I:C^{2n-2}\longrightarrow C^{|I|}$ the one from $C^{2n-2}$.

     We will establish the vanishing of both cycles simultaneously.
    Let us remark that one could also prove the first vanishing directly, viewing $B^n(C, e)$ as a correspondence between $C^2$ and $C^{n - 2}$. Indeed, unwinding the definitions, we have the following equality
    \begin{equation}\label{eq-bnn-2factors}
    B^n(C,e)_*(\Delta_C) = \overline{\pr}_{[3,n],*}\big( B^n(C, e)\cdot \overline{\pr}^*_{1 2}(\Delta_C) \big),
    \end{equation}
    and non-trivial computations (similar to those carried out in this proof) show that $B^n(C, e)_*(\Delta_C)= -\gamma_e^{n-2}(K_C+2e)$.
    Instead, we refine this calculation by computing the cycle
\begin{equation}\label{almostcorresp1}
\overline{\pr}_{[2,n],*}\big( B^n(C,e) \cdot \overline{\pr}^*_{1 2}(\Delta_C) \big).
    \end{equation} 
    The formula for $B^n(C,e)_*(\Delta_C)$ is then deduced by projecting \eqref{almostcorresp1} onto the last $(n-2)$ factors.

    By definition, the cycle \eqref{almostcorresp1} is given by
\begin{equation*}
    \overline{\pr}_{[2, n], *} \big( {\pr_{[n+1,2n],*}} 
        \big( \pr_{[1,n]}^* (\Delta_C^{(n)}) \cdot \pi_1^{\otimes n} \big) \cdot \overline{\pr}_{12}^* (\Delta_C)\big).
\end{equation*}
    Let us write $\pr_{[1,n]}^*(\Delta_C^{(n)})$ as $\pr_{1 2}^* (\Delta_C) \cdot \pr_{[2, n]}^* (\Delta_C^{(n - 1)})$. {Writing $\pr_{[n+2,2n]} = \overline{\pr}_{[2,n]}\circ \pr_{[n+1,2n]}$ and $\pr_{n+1\, n+2} = \overline{\pr}_{12} \circ \pr_{[n+1,2n]}  $, and} applying the projection formula with respect to $\pr_{[n+1,2n]}$, we find that \eqref{almostcorresp1} is equal to
\begin{equation*}
    \pr_{[n+2, 2n], *} \big( \pr_{12}^*(\Delta_C) \cdot \pr_{[2,n]}^*(\Delta_C^{(n-1)}) \cdot \pi_1^{\otimes n} \cdot \pr_{n+1\, n+2}^*(\Delta_C) \big).
\end{equation*}
We replace $\pi_1^{\otimes n}$ by its expression as $$\pr_{\{1,2,n+1,n+2\}}^*(\pi_1^{\otimes 2}) \cdot \pr_{[3,n]\sqcup[n+3,2n]}^*(\pi_{1}^{\otimes (n-2)}).$$
By the definition of $\pi_1$ together with the adjunction formula \eqref{adjunction}, a straightforward computation shows
\[
\pi_1^{\otimes 2} \cdot (\Delta_C\times \Delta_C) = (\Delta_{C}\times \Delta_{C})_* \big( \Delta_{C, *}(-K_C-4e)  + 2 e\times e\big).
\]
Moreover, using the push-pull formula (see \cite[Section 5.9]{MoonenAWS} or \cite[\href{https://stacks.math.columbia.edu/tag/02RF}{Tag 02RF}]{stacks-project}), we have the equality
\[
\pr_{\{1,2,n+1,n+2\}}^*(\Delta_C\times\Delta_C)_* = (\Delta_C\times \id_{C^{n-2}} \times \Delta_C\times \id_{C^{n-2}})_*\widetilde{\pr}_{1n}^*.
\]
Putting everything together, and then applying the projection formula with respect to the map $(\Delta_C\times \id_{C^{n-2}} \times \Delta_C\times \id_{C^{n-2}})\colon  C^{2n-2}\to C^{2n}$,   the cycle~\eqref{almostcorresp1} becomes
\begin{align*}
& \widetilde{\pr}_{[n,2n-2], *} \left( \widetilde{\pr}_{1n}^*\big( \Delta_{C, *}(-K_C-4e)  + 2 e\times e \big)  \cdot \widetilde{\pr}_{[1, n+1]}^* (\Delta_C^{(n+1)}) \cdot {\widetilde{\pr}_{[2, n-1]\sqcup [n+1, 2n-2]}^*}(\pi_1^{\otimes (n-2)}) \right)\\ 
& = \widetilde{\pr}_{[n,2n-2],*} 
\Big( 
\widetilde{\pr}_{[1,n]}^* 
\big( \Delta_{C,*}^{(n)}(-K_C-4e) + \Delta_{C,*}^{(n-1)}(2e)\times e\big)
\cdot 
\widetilde{\pr}_{[2,n-1]\sqcup [n+1,2n-2]}^* (\pi_1^{\otimes (n-2)}) \Big). 
\end{align*}

On the other hand, we note that by definition, the cycle $\beta_e^{n-1}(-K_C-4e)$ is equal to
\begin{align*}
&
\widetilde{\pr}_{[n,2n-2],*} 
\left( 
\widetilde{\pr}_{[1,n-1]}^*\Delta_{C,*}^{(n-1)}(-K_C-4e) 
\cdot 
\widetilde{\pr}_{1 n}^*(\pi_1) 
\cdot 
\widetilde{\pr}_{[2,n-1]\sqcup [n+1,2n-2]}^* (\pi_1^{\otimes (n-2)}) 
\right)
\\
&=
\widetilde{\pr}_{[n,2n-2],*} 
\left( 
\widetilde{\pr}_{[1,n]}^*
\big(
\Delta_{C,*}^{(n)}(-K_C-4e) + \Delta_{C,*}^{(n-1)}(K_C+4e)\times e 
\big)   
\cdot  
\widetilde{\pr}_{[2,n-1]\sqcup [n+1,2n-2]}^* (\pi_1^{\otimes (n-2)}) 
\right). 
\end{align*}
Therefore, the difference $\beta_e^{n-1}(-K_C-4e)-\overline{\pr}_{[2, n], *} (B^n(C, e) \cdot \overline{\pr}_{1 2}^* (\Delta_C))$ is equal to
\begin{align*}
& \widetilde{\pr}_{[n,2n-2],*} 
\Big( 
\widetilde{\pr}_{[1,n-1]}^*\Delta_{C,*}^{(n-1)}(K_C+2e)
\cdot 
\widetilde{\pr}_n^*(e)  
\cdot  
\widetilde{\pr}_{[2,n-1]\sqcup [n+1,2n-2]}^* (\pi_1^{\otimes (n-2)}) 
\Big) \\ &  = e \times \beta_e^{n-2}(K_C+2e).
\end{align*}
In other words, we have found the following simplified expression for \eqref{almostcorresp1}:
\begin{equation}\label{eq-33simplified}
    \overline{\pr}_{[2,n],*}\big( B^n(C,e) \cdot \overline{\pr}^*_{1 2}(\Delta_C) \big)  = \beta_e^{n-1}(-K_C-4e)-(e \times \beta_e^{n-2}(K_C+2e)).
\end{equation} 

By the expression for $\beta_e^{n-1}$ given in \Cref{explicitexpression}, it is straightforward to see that the projection onto the last $(n-2)$ factors of $\beta_e^{n-1}(-K_C-4e)$ is trivial. Therefore by \eqref{eq-bnn-2factors} we find that $$B^n(C,e)_*(\Delta_C)=-\beta_e^{n-2}(K_C+2e).$$
Combining this with \eqref{eq-33simplified}, we deduce that the vanishing of $B^n(C,e)$ implies the vanishing of both $\beta^{n-1}_e(K_C+4e)$ and $\beta_e^{n-2}(K_C+2e)$. The statement of the proposition then follows after recalling that $\beta_e^{m}$ and $\gamma_e^{m}$ agree on zero-cycles. 
\end{proof}

Viewing $B^n(C,e)$ as a correspondence from $C$ to $C^{n-1}$, we deduce another consequence of $B^n(C,e)=0$.

\begin{proposition}\label{vanish3}
Let $n \geq 3$. Then for any $z\in \CH_0(C)$, we have
\[
B^{n}(C,e)_*(z) = \gamma_e^{n-1}(z) - \deg(z)\cdot \gamma_e^{n-1}(e).
\]
In particular by \Cref{vanishing of gamma n-12}, if $B^n(C,e)=0$, then $\gamma_e^{n-1}(z)=0$.
\end{proposition}

\begin{proof}
{As in the previous proposition we denote by $\pr_I$ the projections from $C^{2n}$, by $\overline{\pr}_I$ the projections from $C^{n}$, and by $\widetilde{\pr}_I$ the projections from $C^{2n-2}$.}

We view $B^n(C, e)$ as a correspondence from $C$ to $C^{n - 1}$, and compute $B^n(C, e)_* (z)$.
    As before, this can be written as $\pr_{[n + 2, 2n], *}\big(\pi_1^{\otimes n}\cdot  \pr_{n + 1}^* (z)\cdot \pr_{[1, n]}^* (\Delta_C^{(n)}) \big)$.
    Using that $\pi_1=\Delta_C - (e \times C) - (C \times e)$, we obtain 
\begin{equation*}
    B^n(C, e)_*(z) = \pr_{[n + 2, 2n], *}
    \big(\pr_{[2, n] \sqcup [n + 2, 2n]}^*(\pi_1^{\otimes (n-1)}) \cdot  
    \pr_{1\, n + 1}^* \big(\Delta_{C,*}(z) - e \times z
    \big)     \cdot 
    \pr_{[1, n]}^* (\Delta_C^{(n)})
    \big).
\end{equation*}
We then can write $\pr_{1\, n + 1}^* (\Delta_{C,*}(z) - e \times z) \cdot \pr_{[1, n]}^* (\Delta_C^{(n)})$ as 
\begin{equation*}
     \pr_{1\, n + 1}^* (\Delta_C)\cdot \pr_{[1, n]}^* (\Delta_{C,*}^{(n)}(z)) - \pr_{n + 1}^* \left(z\right)\cdot \pr_{[1, n]}^* (\Delta_{C,*}^{(n)}(e))
\end{equation*}
    and obtain that $B^n(C, e)_*(z)$ is equal to 
    \begin{align*}
    \pr_{[n + 2, 2n], *}\big(\pr_{[2, n] \sqcup [n + 2, 2n]}^*(\pi_1^{\otimes (n - 1)})\cdot\pr_{1\, n + 1}^*(\Delta_C) \cdot  &\pr_{[1, n]}^* (\Delta_{C,*}^{(n)}(z))\big) \\- \pr_{[n + 2, 2n], *}\big(\pr_{[2, n] \sqcup [n + 2, 2n]}^*(\pi_1^{\otimes (n - 1)}) &\cdot \pr_{n + 1}^*(z) \cdot \pr_{[1, n]}^* (\Delta_{C,*}^{(n)}(e))\big).
    \end{align*}
    We factor the the diagonal $\Delta^{(n)}_{C,*}(z)$ as $\overline{\pr}_{12}^*(\Delta_C) \cdot \overline{\pr}_{[2,n]}^*(\Delta_{C,*}^{(n-1)}(z))$, and write $\pr_{[n + 2, 2n]}$ as $\widetilde{\pr}_{[n, 2n-2]}\circ \pr_{[2, n] \sqcup [n + 2, 2n]}$. Applying the projection formula with respect to $\pr_{[2, n] \sqcup [n + 2, 2n]}$, we see that $B^n(C, e)_*(z)$ is equal to
    \begin{align*}
    \widetilde{\pr}_{[n, 2n-2], *}
    \big( \widetilde{\pr}_{[1, n-1]}^* (\Delta_{C,*}^{(n-1)}(z)) \cdot \pi_1^{\otimes (n - 1)} \big) 
    - 
    \widetilde{\pr}_{[n , 2n-2], *}
    \big(
    \deg(z) \cdot  \widetilde{\pr}_{[1, n-1]}^* (\Delta_{C,*}^{(n-1)}(e))\cdot \pi_1^{\otimes (n - 1)} 
    \big)
    \\
    = 
    \beta_e^{n-1}(z) - \deg(z)\cdot \beta_e^{n-1}(e) 
    = \gamma_e^{n-1}(z) - \deg(z)\cdot \gamma_e^{n-1}(e),
    \end{align*}
    as desired. \qedhere
\end{proof}

We are now in a position to prove the main result of this section.

\begin{theorem}\label{thm:gammavanishiffB}
For $n \geq 2$, we have $\Gamma^n(C, e)  = 0$ if and only if $B^n(C, e) = 0$.
    \label{BimpliesGamma}
\end{theorem}
\begin{proof}
The forward direction is clear by the remarks at the beginning of this subsection. 
For the backwards direction, if $n=2$, then $B^2(C,e)$ is equal to $\Gamma^2(C,e)$ as noted in \Cref{examples}. 
If $B^n(C, e) = 0$ and $n\geq 3$, then the second {part} of \Cref{vanish3} implies that $\gamma_e^{n-1}(e)=0$.
The difference $\Gamma^n(C, e) - B^n(C, e)$ is precisely $\sum_{i = 1}^{n} \widehat{\pr_{i}}^* \gamma_e^{n - 1}(e)$, and so we have $\Gamma^n(C,e)=0$.
\end{proof}

\begin{remark}
    {For $n=3$, the fact that the vanishing of $\Gamma^3(C,e)$ implies the vanishing of $\gamma_e^2(e)$ has also been observed by Qiu and Zhang \cite[Proposition 2.3.2]{QZ24}.}
\end{remark}

\subsection{Successive Vanishing}\label{sec-successive vanishing}
By \cite[Prop. 3.4]{MY16}, if $e=x_0$ is a $k$-point and $z\in \CH(C)$ any cycle, then the vanishing of $\gamma_{x_0}^n(z)$ forces the vanishing of $\gamma_{x_0}^{n+1}(z)$. However the method in \emph{loc. cit.} does not extend to the case of a general divisor $e$.

In this section we {prove a similar result under a slightly stronger assumption: for any divisor $e$ of degree 1, the vanishing of $\gamma_e^n(z)$ and $\gamma_e^{n+1}(z)$ implies the vanishing of all the following $\gamma_e^k(z)$ for $k\geq n+2$.}
To do this, we establish a relationship between $\gamma_e^{n + m}(z)$ and $\gamma_e^{n - 1}(z)$, under the assumption $\gamma_e^n(z)=0$.

The main idea is to replace $\Delta_C^{(n+m)}(z)$ with {$\pr_{[1, n]}^* \Delta_C^{(n)}(z) \cdot \pr_{[n, n + m]}^* \Delta_C^{(m + 1)}$}, and then use the vanishing of $\gamma_e^n(z)$ to rewrite the first term. This could be done using the formulae in \Cref{explicitexpression}, but as in the previous section, this is notationally complicated, and the proof is simplified using projectors.
\begin{theorem}\label{thm:recursiverelationbetagamma}
    Let $z\in \CH(C)$ be such that $\gamma_e^n(z) = 0$,
    and let $n\geq 2$.
    Then for all $m \geq 1$, $$\gamma_e^{n + m}(z) = \gamma_e^{n - 1}(z) \times \gamma_e^{m + 1}(e).$$ Similarly, if $\beta_e^n(z) = 0$, then $$\beta_e^{n + m}(z) = \beta_e^{n - 1}(z) \times \beta_e^{m + 1}(e) + \beta_e^{n - 1}(z \cdot e) \times \beta_e^{m + 1}(C).$$
\end{theorem}
\begin{proof}
By definition we have $\gamma_e^{n + m}(z) = \pi_{+, *}^{\otimes (n + m)} \Delta^{(n + m)}_{C, *}(z)$. This expression can be written as
\begin{equation}\label{firsteq}
\pi_{+, *}^{\otimes (n + m)} \big(
\pr_{[1,n]}^*\Delta^{(n)}_{C, *}(z) \cdot \pr_{[n,n+m]}^*\Delta^{(m+1)}_C  \big).
\end{equation}
We will first rewrite $\Delta^{(n)}_{C, *}(z)$ as a sum of explicit terms.
Then we will plug in this ``new expression'' of $\Delta_{C,*}^{(n)}(z)$ in equation~\eqref{firsteq} and see that only one term contributes in a non-trivial way. 
We then finish by interpreting that term as $\gamma_e^{n-1}(z)\times \gamma_e^{m+1}(e)$.

Since $\gamma_e^n(z) = 0$, by the definition of $\gamma_e^n(z)$ we have  
\[
\Delta^{(n)}_{C, *}(z) = \Delta^{(n)}_{C, *}(z) - \gamma_e^{n}(z) = (\Delta_C^{\otimes n} - \pi_+^{\otimes n})_* \big(\Delta^{(n)}_{C, *}(z)\big) . 
\]
Writing $\pi_+$ as $\Delta_C - (C \times e)$, we can expand the difference $\Delta_C^{\otimes n} - \pi_+^{\otimes n}$ as

\[
\sum_{I \subsetneq \{1,\dots,n\}} (-1)^{n-1-|I|}\cdot \prod_{i\in I} \pr_{i\, n+i}^*(\Delta_C) \cdot \prod_{i\notin I} \pr_{i\, n+i}^*(C\times e).
\]
Therefore, {the pushforward $(\Delta_C^{\otimes n} - \pi_+^{\otimes n})_* (\Delta^{(n)}_{C, *}(z) )$ coincides with the sum over $I \subsetneq \{1,\dots,n\}$ of the terms 
\begin{equation}\label{eq:terminthepushedfforward}
(-1)^{n-1-|I|} \cdot \big(\prod_{i\in I} \pr_{i\, n+i}^*(\Delta_C) \cdot \prod_{i\notin I} \pr_{i\, n+i}^*(C\times e)\big)_* \big(\Delta_{C, *}^{(n)}(z)\big).
\end{equation} 
Carrying on the computations,  one can see that equation~\eqref{eq:terminthepushedfforward} 
is equal to
\begin{equation}\label{eq:casesInonemptyorempty}
\begin{cases}
    (-1)^{n-1}\deg(z)\cdot \prod_{i=1}^n \overline{\pr}_i^*(e) & \textup{if } I=\varnothing,
    \\
    (-1)^{n-1-|I|}\cdot \overline{\pr}^*_{I}(\Delta_{C,*}^{(|I|)})(z) \cdot \prod_{i\notin I} \overline{\pr}_i^*(e) & \textup{if } I\neq \varnothing,
\end{cases} 
\end{equation}
where $\overline{\pr}\colon C^n\to C$ is the projection onto the $i$-th component. Since $I\subsetneq \{1,\ldots,n\}$ is a proper subset, this expression always has the shape $\widehat{\overline{\pr}}_j^*(z_j')\cdot \overline{\pr}_j^*(e)$ for some $ j \notin I$ and  cycle $z_j'\in \CH(C^{n-1})$.}
We now note that for any $j\neq n$, the pushforward $\pi_{+, *}^{\otimes (n + m)} (\widehat{\pr}_j^*(z_j')\cdot {\pr}_j^*(e) \cdot {\pr_{[n,n+m]}^*\Delta^{(m+1)}_C})$ is zero, as $\pi_{+, *}(e) = 0$. 
Therefore the only non-trivial contribution comes from $I=\{1,\dots,n-1\}$.
{
Hence, the original equation \eqref{firsteq} is equal to $\pi_{+, *}^{\otimes (n+m)} ( \Delta_{C,*}^{(n-1)}(z) \times \Delta_{C,*}^{(m+1)}(e) )$, which is exactly $\gamma_e^{n - 1}(z) \times \gamma_e^{m + 1}(e)$.}

The argument for $\beta_e$ is analogous. We use $\pi_1 = \pi_+ - \pi_2$ and the assumption $\beta^{n}_e(z)=0$ to rewrite   \[\Delta_{C,*}^{(n)}(z) =  \left(\Delta_C^{\otimes n} - (\pi_+ - \pi_2)^{\otimes n}\right)_*\Delta_{C,*}^{(n)}(z).\] For dimension reasons, $\pi_{2, *}^{\otimes j}(\Delta_{C, *}^{(j)}(z'))=0$ for any $z'\in \CH(C)$ and $j\geq 2$. Hence the only non-zero terms in this expression come from pushforward by the cycle
\begin{equation*}
\Delta_C^{\otimes n}-\big( \pi_+^{\otimes n} - \sum_{i=1}^n \widehat{\pr}^*_{i\, n+i}(\pi_+^{\otimes (n-1)})\cdot \pr^*_{i\, n+i}(\pi_2)\big).
\end{equation*}
Therefore, the cycle $\Delta_{C, *}^{(n)}(z)$ is equal to
\begin{align*}
 (\Delta_C^{\otimes n} - \pi_+^{\otimes n})_* \Delta_{C, *}^{(n)}(z)
 + \sum_{i=1}^n \pr_{[n+1,2n], *} \big( \pr_{[1,n]}^*\big(\Delta_{C, *}^{(n)}(z\cdot e) \big) \cdot \widehat{\pr}^*_{i\, n+i}(\pi_+^{\otimes (n-1)})\cdot \pr_{n+i}(C) \big).
\end{align*}
As done before, we then plug in this new expression of $\Delta_{C, *}^{(n)}(z)$ in the equation computing $\beta_e^{n+m}(z)$ (that is, \eqref{firsteq} with $\pi_+$ replaced by $\pi_1$). The first summand contributes in the same way as 
$\Delta_{C, *}^{(n-1)}(z)\times e$; the only non-trivial contribution from the second sum is from the term $i=n$,  because $\pi_{1, *}(C)=0$.
Therefore this second piece contributes in the same way as $\Delta_{C, *}^{(n-1)}(z\cdot e)\times C$. \qedhere
\end{proof}

As an immediate application, we get the first statement of \Cref{maintheoremB}.
\begin{corollary}\label{successivegammavanishing}
    Let $n\geq 1$. 
    If $z \in \CH(C)$ is such that both $\gamma_e^{n}(z)$ and $\gamma_e^{n+1}(z)$ vanish, then $\gamma_e^{k}(z) = 0$ for all $k \geq n$.
    \label{seqvanishgamma}
    In particular for $z=[C]$, we get 
    \[
    \Gamma^{n}(C,e)= \Gamma^{n+1}(C,e) = 0\quad  \implies\quad  \Gamma^k(C,e)=0\quad \forall k \geq n.
    \]
\end{corollary}

We also have the following similar result under weaker assumptions.

\begin{corollary}\label{corollaryonvanishingofgamma}
    Let $n\geq 2$.
    If $\gamma_e^{n}(e) = 0$, then $\gamma_e^{k(n - 1) + 1}(e) = 0$ for all $k \geq 0$.
\end{corollary}

\begin{proof}
We proceed by induction. The case $k=0$ is always true,  and the case $k=1$ is the hypothesis.
Let us then assume that $\gamma_e^{k(n-1)+1}(e)=0$ for $k\geq 1$. By \Cref{thm:recursiverelationbetagamma}, we have
\begin{equation*}
    \gamma_e^{k(n-1)+1+m}(e) = \gamma_e^{k(n-1)}(e)\times \gamma_e^{m+1}(e), \quad \forall \ m\in \mathbb{N}.
\end{equation*}
Setting $m=n-1$ we get $\gamma_e^{(k+1)(n-1)+1}(e)=0$.
\end{proof}

\begin{theorem}\label{thm:vanishinglittlegamma}
    Assume $\Gamma^n(C, e) = 0$ for some $n \geq 3$. Then for any $z \in \CH_0(C)$, we have $\gamma^k_e(z) = 0$ for all $k \geq n - 1$.
\end{theorem}
\begin{proof}
The case when $k = n - 1$ is given by \Cref{vanish3} and \Cref{thm:gammavanishiffB}. By \Cref{successivegammavanishing}, it is sufficient to prove this for $k = n$.

Let us denote by $\pr_I$ the $I$-projection from $C^n$. 
As before, we write $\gamma_e^n(z) = \pi_{+, *}^{\otimes n} (\Delta^{(n)}_{C,*}(z))$. We now rewrite $\Delta^{(n)}_{C,*}(z)$ as $\pr_1^* (z) \cdot \Delta^{(n)}_{C}$. Since $\Gamma^n(C, e) = 0$, we have $\pi_{+, *}^{\otimes n} \Delta^{(n)}_C = 0$, and so $\Delta^{(n)}_C = (\Delta^{\otimes n}_C - \pi_+^{\otimes n})_* \Delta^{(n)}_C$. 
The expression \eqref{eq:casesInonemptyorempty} from the proof of \Cref{thm:recursiverelationbetagamma} tells us that every term in this pushforward can be written as a product of the form $\widehat{{\pr}}_j^{*}(z'_j)\cdot {\pr}_j^*(e)$ for some cycle $z'_j\in \CH(C^{n-1})$ and some index $1 \leq j \leq n$. 
Intersecting this with $\pr_1^* (z)$ does not change this fact. 
Indeed either a term contains ${\pr}_j^*(e)$ for $j\neq 1$; or it contains ${\pr}_1^*(e)$, and the intersection then vanishes for dimension reasons. 
By \eqref{decompositionofpushforward}, the remaining terms will vanish under $\pi_{+, *}^{\otimes n}$ as $\pi_{+, *}(e) = 0$.
\end{proof}

We can leverage this to get a version of \Cref{corollaryonvanishingofgamma}, for $\Gamma^n(C, e)$. This gives the second statement of \Cref{maintheoremB}.
\begin{corollary}
    Let $n \geq 3$, and suppose $\Gamma^n(C, e) = 0$. Then for all $k \geq 2n - 2$, $\Gamma^{k}(C, e) = 0$.
    \label{vanishinglargegamma}
\end{corollary}
\begin{proof}
    Let $k = n + m$, where $m \geq n - 2$. By \Cref{thm:recursiverelationbetagamma}, $\Gamma^{k}(C, e) = \Gamma^{n - 1}(C, e) \times \gamma_e^{m + 1}(e)$. By the previous result, $\gamma_e^{m + 1}(e) = 0$. 
\end{proof}

\begin{remark}\label{altrecurrence}
    The relations written above are not the only possible recurrence relations. 
    For example, we can write $\Delta^{(n + 2)}_C$ as $\pr_{1 2}^* \Delta_C \cdot \pr_{[2, n + 1]}^* \Delta^{(n)}_C  \cdot \pr_{n + 1\, n + 2}^* \Delta_C$.
    If for a $0$-cycle $z$ we have $\gamma_e^n(z) = 0$, then with similar steps as for \Cref{thm:recursiverelationbetagamma} we would get
    \[
        \gamma_e^{n + 2}(z) = - \gamma_e^2(e) \times \gamma_e^{n - 2}(z) \times \gamma_e^2(e) \quad \forall n\geq 3.
    \]
    Using similar arguments, we could get additional consequences of the vanishing of $\Gamma^n(C, e)$ (or equivalently, $B^n(C, e)$).
\end{remark}

\section{From Jacobians to Products of Curves and Back Again}\label{JandproductofC} 
The aim of this section is to prove \Cref{maintheoremA}. We recall that this theorem is proven in \cite{MY16} in the case of a base point, but their method does not generalize to arbitrary $e$ due to the fact that we do not necessarily have the vanishing of $\gamma_e^2(e)$ (or equivalently $\updelta^e_{(2)}$, see \Cref{equivalentvanishingforzerocycles}).

In the course of the proof, we will see that the cycle $B^{s+2}(C,e)$ naturally appears in terms of the Beauville component  $[C]_{(s)}$. 
This could be anticipated by the fact that the motive $\mathfrak{h}_1(C)$ corresponds to the Jacobian (see \cite[\S 3.2]{Scholl94}), and so the K\"unneth decomposition for motives suggests that cycles on $C^n$ related to the Jacobian should occur in $\mathfrak{h}_{1 \ldots 1}(C^n)$. This is one reason for introducing this cycle.

As a consequence, we show that the vanishing of $[C]_{(s)}$ implies the vanishing of $B^{s+2}(C,e)$ (see \Cref{ceresatomoddiag}), which in turn is equivalent to the vanishing of $\Gamma^{s+2}(C,e)$ by \Cref{thm:gammavanishiffB}.
We then show that the vanishing of $\Gamma^{s+2}(C,e)$ implies the vanishing of $[C]_{(s)}$ in \Cref{moddiagtoceresa}, which yields \Cref{maintheoremA}. 

\begin{notation}\label{notationsection4}
In this section, we label the coordinates in $C\times C^n$ by $\{0,1,\dots ,n\}$ and the projections from it by $\overline{\pr}$. Similarly, we label the coordinates in $J \times J^n$ by $\{0,1,\ldots,n\}$ and the projections from it by $\overline{\pr}^J$.
We abuse the notation and denote the projections from $C^n$
(so forgetting first the $0$-th coordinate) also by $\overline{\pr}$, and similarly for $\overline{\pr}^J$. As in the previous section, we reserve the symbol $\pr$ for projections from $C^{2n}$.
\end{notation}
\subsection{From the Beauville Decomposition to Diagonal Cycles}
In this subsection we show that for any $z \in \CH_i(C)$ and $s \geq 1$, the vanishing of $\iota_{e, *} (z)_{(s)}$ in $\CH_{i, (s)}(J)$ implies $\beta_e^{s + 2i}(z) = 0$ in $\CH(C^{s + 2i})$. To achieve this, we follow and extend the proof method of \cite[Theorem 5.3.1]{Zhang10}. Specifically, our goal will be to compute $(\phi \circ \sigma_{s + 2i})^* \mathscr{F} (\iota_{e, *} (z)_{(s)})$, where $\sigma_{n}$ denotes the composite 
$$\sigma_n\colon C^n \overset{\iota_e^n}{\longrightarrow} J^n \xrightarrow[]{m_{[1,n]}} J.$$ 
Let us recall that $\sigma_n$ can also be described as 
\[
C^n \overset{a_n}{\longrightarrow} \Pic^n_C  \xrightarrow[]{t_{-ne}} \Pic^0_C=J,
\]
where $a_n$ and $t_{-ne}$ are as in \Cref{def:translation} and the discussion immediately preceding it.

\begin{proposition}\label{expression of fourier transform} 
Let $z\in \CH(C)$.
    We have the following equality in  $\CH(C^n)$
    $$
    (\phi \circ \sigma_{n})^* \mathscr{F} (\iota_{e, *} (z)) = 
    \overline{\pr}_{[1,n], *} \big( \overline{\pr}_0^* z \cdot \ch\big((\iota_e \times \sigma_{n})^* c_1(\Lambda_\Theta)\big) \big). $$
  
\end{proposition}
\begin{proof}
    We have the following Cartesian diagram:
    $$\begin{tikzcd}[column sep={14mm}, row sep={10mm}]
        C 
        \arrow[d, swap, "{\iota_e}"] 
        & 
        C \times J^t 
        \arrow[l, swap, "{\pr_C}"] 
        \arrow[d, swap, "{\iota_e \times \id_{J^t}}"] 
        & C \times J 
        \arrow[l, swap,  "{\id_C \times \phi}"] 
        \arrow[d, "{\iota_e \times \id_J}"'] 
        & C \times C^n 
        \arrow[l, swap, "{\id_C \times \sigma_{n}}"] 
        \arrow[d, "{\iota_e \times \id_{C^n}}"'] 
        \\
        J 
        & 
        J \times J^{t} 
        \arrow[l, "{\pr_J}"'] 
        \arrow[d, swap, "{\pr_{J^t}}"] 
        & J \times J 
        \arrow[l, swap, "{\id_J \times \phi}"] 
        \arrow[d, swap] 
        & J \times C^n 
        \arrow[l, "{\id_J \times \sigma_n}"'] 
        \arrow[d, "{\pr_{C^n}}"']
        \\
        & 
        J^t 
        & 
        J 
        \arrow[l, "{\phi}"'] 
        & 
        C^n .
        \arrow[l, "{\sigma_n}"']
    \end{tikzcd}$$
    We start with the maps at the bottom left of the diagram. Since all the squares are Cartesian, the {push-pull formula} allows us to work our way to the top right: 
    \begin{align*}
        \sigma_{n}^* \phi^*
        \left( \pr_{J^t,*}\big(\pr_J^*(\iota_{e, *} (z)) \cdot \ch(\mathscr{P}) \big)
        \right) &
        = \pr_{C^n, *} \left( (\id_J\times \phi\circ\sigma_n)^* \big(\pr_J^*(\iota_{e, *} (z)) \cdot \ch(\mathscr{P}) \big) \right)
        \\
        & = \pr_{C^n, *} \big( (\iota_e\times \id_{C^n})_*(z\times C^n) \cdot (\id_J \times \sigma_n)^* \ch( \Lambda_{\Theta} ) \big)
        \\ 
        & = \overline{\pr}_{[1,n], *} \left( \overline{\pr}_0^* z\cdot \ch\big( (\iota_e \times \sigma_{n})^*\Lambda_{\Theta} \big) \right).
    \end{align*}
    In the last two steps, we first use that $\Lambda_{\Theta}=(\id_J\times \phi)^*\mathscr{P}$, and then apply the projection formula with respect to $\iota_e\times \id_{C^n}$. 
    This yields the desired equality.
\end{proof}
Now we would like to rewrite $(\iota_e \times \sigma_{n})^*\Lambda_{\Theta}$ in other terms.
\begin{proposition}\label{pullbackofmumford} 
    We have the following equality in $\CH(J\times C^n)$,
    \[
   (\id_J \times \sigma_{n})^* c_1(\Lambda_{\Theta}) = 
    \sum_{j = 1}^{n} 
    (\id_J\times \overline{\pr}_j)^* (\id_J\times \iota_e)^*c_1(\Lambda_{\Theta}),
    \]
    and therefore the following one in $\CH(C\times C^n)$,
    \[
    (\iota_e \times \sigma_{n})^* c_1(\Lambda_{\Theta}) = -\sum_{j = 1}^{n} \overline{\pr}_{0 j}^* \pi_1.
    \]
\end{proposition}

\begin{proof}
    Let $D = m_{[0, n]}^* \Theta_{\kappa} - m_0^* \Theta_{\kappa} - m_{[1, n]}^* \Theta_{\kappa}$, then we have the following equality in $\CH(J\times C^n)$
    \begin{align*}
    (\id_J\times \iota_e^{n})^*
    D &= 
    (\id_J\times \iota_e^{n})^*m_{[0,n]}^*
    \Theta_{\kappa} - 
    \pr_J^*
    \Theta_{\kappa} - \pr_{C^n}^* \sigma_n^*
    \Theta_{\kappa} 
    \\
    &= 
    (\id_J\times \sigma_{n})^*m^*
    \Theta_{\kappa} - 
    \pr_J^*
    \Theta_{\kappa} - \pr_{C^n}^* \sigma_n^*
    \Theta_{\kappa} 
    \\ &
    = 
    (\id_J \times \sigma_n)^* (m^* - \pr_{J, 1}^* - \pr_{J, 2}^*) \Theta_{\kappa} 
    =(\id_J \times \sigma_n)^*c_1(\Lambda_{\Theta}).
    \end{align*}
    Observe that $m_{[0,k]}$ can be written as the sum of $m_0$, $m_{[1,k - 1]}$ and $m_k$.
    Applying the Theorem of the Cube to the triple $(m_0, m_{[1,k - 1]}, m_k)$ then gives 
\begin{equation*}
            (m_{[0, k]}^* - m_{[1, k]}^*) \Theta_{\kappa}  = (m_{[0, k - 1]}^* - m_{[1, k - 1]}^*) \Theta_{\kappa} + (m_{0 k}^* - m_0^* - m_k^*)\Theta_{\kappa}.
\end{equation*}
    By induction on $k$, we can write 
\begin{equation*}
            (m_{[0, k]}^* - m_{[1, k]}^*) \Theta_{\kappa}  = (m_{0 1}^* - m_{1}^*) \Theta_{\kappa} + \big( \sum_{j=2}^k m_{0 j}^*  - (k-1)m_0^* - \sum_{j=2}^k m_j^* \big)\Theta_{\kappa}.
\end{equation*}
    Since $D+m_0^*\Theta_{\kappa}=(m_{[0, n]}^* - m_{[1, n]}^*) \Theta_{\kappa}$, we get
 \begin{equation*}
      D = \big( \sum_{j=1}^n m_{0 i}^*  - n \, \overline{\pr}_{ 0}^{J,*} - \sum_{j=1}^n \overline{\pr}_{j}^{J,*} \big) \Theta_{\kappa} = { 
        \sum_{j = 1}^n \overline{\pr}_{0j}^{J,*} (c_1(\Lambda_{\Theta}))}.
 \end{equation*}
    Therefore, the pullback $(\id_J\times \iota_e^{n})^* D$ can be written as $\sum_{j = 1}^{n} 
    (\id_J\times \overline{\pr}_j)^* (\id_J\times \iota_e)^*c_1(\Lambda_{\Theta})$. Hence the first equality. The second one follows by taking the pullback of the first one by $(\iota_e\times \id_{C^n})$ and using the identity $(\iota_e\times \iota_e)^*c_1(\Lambda_{\Theta})=-\pi_1$ of \Cref{lem-pullbackMB}.
\end{proof}

We can now compute the value of $(\phi \circ \sigma_{s + 2i})^* {\mathscr{F}(\iota_{e, *} (z)_{(s)})}$. We record the following proposition to simplify the intersection calculations.

\begin{lemma}\label{betaz}
Let $z \in \CH(C)$. We have the following equality in $\CH(C^n)$:
\begin{equation*}
      \overline{\pr}_{[1,n], *}\big(\overline{\pr}_0^*\, z \cdot \overset{n}{\underset{j = 1}{\prod}} \overline{\pr}_{0 j}^*\, \pi_1\big) = \beta_e^n(z).
\end{equation*}
\end{lemma}
\begin{proof}     
We first note that $\overline{\pr}_0^*z$ is equal to $\pr_{[n,2n],*} \pr_{[1,n]}^*\Delta_{C,*}^{(n)}(z)$ in $\CH(C\times C^{n})$. Applying the projection formula with respect to $\pr_{[n,2n]}$, the cycle \begin{align*}
    \overline{\pr}_{[1,n], *} \big( \pr_{[n,2n],*} \pr_{[1,n]}^*\Delta_{C,*}^{(n)}(z) \cdot \overset{n}{\underset{j = 1}{\prod}} \overline{\pr}_{0 j}^* \pi_1\big)
\end{align*} can be rewritten as 
\begin{equation*}
\pr_{[n + 1,2n], *}\big( 
\pr_{[1,n]}^*\Delta_{C,*}^{(n)}(z) 
\cdot 
\overset{n}{\underset{j = 1}{\prod}} \pr_{n\, n + j}^*\pi_1\big).
\end{equation*}
Since $\pr_{[1, n]}^* \Delta_C^{(n)} \cdot \pr_{n\, n + j}^* \pi_1 = \pr_{[1, n]}^* \Delta_C^{(n)} \cdot \pr_{j\, n + j}^* \pi_1$, the previous expression is equal to
\begin{equation*}
    \pr_{[n + 1, 2n], *}\big(\pr_{[1, n]}^* \Delta^{(n)}_{C,*}(z) \cdot \prod_{j = 1}^n \pr_{j\, n + j}^* \pi_1\big)
\end{equation*}
which is precisely $\beta_e^n(z)$.
\end{proof}

\begin{theorem}\label{ytogamma}
Let $z\in \CH_0(C)$ and $n\geq 1$.
Then 
\[
(\phi \circ \sigma_{n})^* \mathscr{F}(\iota_{e, *}( z)_{(n)}) = (-1)^n \beta_e^{n}(z)= (-1)^n \gamma_e^{n}(z).
\]
\end{theorem}
\begin{proof}
    By \Cref{expression of fourier transform} and  \Cref{pullbackofmumford}, we have the equality  
    \begin{equation*}
        (\phi \circ \sigma_{n})^* \mathscr{F}(\iota_{e, *} (z))
        = 
        \overline{\pr}_{[1, n], *}
        \big(
        \overline{\pr}_0^* z 
        \cdot 
        \ch\big(-\sum_{j = 1}^n \overline{\pr}_{0 j}^* \pi_1\big) \big), 
    \end{equation*}
    where as before $\ch$ denotes the Chern character. The $n$-th Beauville component of $\iota_{e, *} (z)$ becomes the codimension $n$ component under the Fourier transform. 
    Thus we just have to take the codimension $n$ component of the right hand side, which comes from the codimension $n$ term of $\ch(-\sum_{j = 1}^n \overline{\pr}_{0 j}^*\pi_1)$. 
    Write $\ch(-\sum_{j=1}^n \overline{\pr}_{0 j}^*\pi _1)$ as $(-1)^n \prod_{j=1}^n \ch(\overline{\pr}_{0 j}^*\pi _1)$.
    As $\pi_1$ has codimension one in $C^2$, we have $\pi_1^3 = 0$. 
    Similarly, we have  $\overline{\pr}_0^* z \cdot \overline{\pr}_{0 i}^* \pi_1^2 = 0$. 
    Therefore, we can restrict to the terms of the Chern character with no quadratic or higher terms in any $\overline{\pr}_{0 i}^* \pi_1$, which is precisely the term $(-1)^n \prod_{j = 1}^n \overline{\pr}_{0 j}^* \pi_1$. It only remains to compute $\overline{\pr}_{[1,n], *}\big(\overline{\pr}_0^* z \cdot \prod_{j = 1}^n \overline{\pr}_{0 j}^* \pi_1 \big)$, which is the content of \Cref{betaz}.
\end{proof}

We now fix $n=s+2$ in \Cref{notationsection4}. 

\begin{theorem}\label{pullback of fourier transform}
    The cycle $(\phi \circ \sigma_{s + 2})^* \mathscr{F}([C]^e_{(s)})$ is equal to 
\begin{equation}\label{eq:equalityforFourier(C)}
    (-1)^{s}
    \Big(
    B^{s + 2}(C, e) - \frac{(s + 1)}{2} \sum_{j = 1}^{s + 2} \widehat{{\overline{\pr}}}_j^* \beta_e^{s + 1}(K_C + 4e) -   
    \frac{1}{2} \sum_{j = 1}^{s + 2} \sum_{k \neq j} 
    {\overline{\pr}}_k^* (e) \cdot \widehat{{\overline{\pr}}}_{j k}^* \beta_e^{s}(K_C + 2e)
    \Big).
    \end{equation}
\end{theorem}
\begin{proof}
Set $z=[C]$ and let $s\ge 0$. By \Cref{expression of fourier transform}, we have the equality 
\begin{equation*}
    (\phi \circ \sigma_{s + 2})^* \mathscr{F}([C]^e) 
= 
\overline{\pr}_{[1,s + 2], *}
\big(
\ch\big(
-\sum_{i = 1}^{s + 2} \overline{\pr}_{0i}^* \pi_1
\big)
\big).
\end{equation*}
The $s$-th Beauville component of $[C]^e$ becomes the codimension $(s + 1)$ component under the Fourier transform. Hence, we only need to compute the codimension $(s + 1)$ component of the right hand side, which comes from the codimension $(s + 2)$ 
component of $\ch(-\sum_{i = 1}^{s + 2} \overline{\pr}_{0  i}^*\pi_1)$. As before, since $\pi_1$ has codimension one in $C^2$, we have $\pi_1^3 = 0$. 
Unlike in the previous case, we do get contributions from $\overline{\pr}^{*}_{0 i} (\pi_1^2)$; however we can use the fact that $\overline{\pr}^{*}_{0 i} (\pi_1^2) \cdot \overline{\pr}^{*}_{0 j} (\pi_1^2) = 0$ for dimension reasons.  
Therefore, we can restrict to the terms of the Chern character with at most one quadratic term in any $\overline{\pr}_{0 j}^* \pi_1$. 
These are given by the expression
\begin{equation}\label{doublesum}
    (-1)^{s+2}
\prod_{j = 1}^{s + 2} \overline{\pr}_{0 j}^* \pi_1 + \frac{(-1)^{s}}{2} \sum_{j = 1}^{s + 2} \sum_{k \neq j} \big(\overline{\pr}_{0 k}^* \pi_1^2 \prod_{l \neq j, k} \overline{\pr}_{0 l}^* \pi_1\big).
\end{equation}
The first term contributes $(-1)^{s+2}\beta_e^{s + 2}(C) = (-1)^{s+2} B^{s + 2}(C, e)$ by \Cref{betaz}. 
Hence it only remains to compute the double sum. 

By the adjunction formula \eqref{adjunction} and the definition of the cycle $\pi_1$, we find that $$\pi_1^2= - \Delta_{C,*}(K_C +4e) + 2\cdot(e\times e).$$
Writing $\Delta_{C,*}(K_C +4e)$ as $\Delta_C\cdot((K_C+4e)\times C)$, and in turn $\Delta_C$ as $\pi_1+(C\times e)+(e\times C)$, we have that 
$$ \pi_1^2=-\pi_1\cdot \big((K_C+4e)\times C\big)-(K_C+2e)\times e. $$
Using this expression, we find that the double sum in \eqref{doublesum} is equal to
\begin{multline*}
\frac{(-1)^{s+1}}{2} \sum_{j = 1}^{s + 2} \sum_{k \neq j} 
\Big[
\overline{\pr}_{[1,s + 2], *}
\big(\overline{\pr}_{0}^* (K_C + 4e) \cdot \prod_{l \neq j} \overline{\pr}_{0 l}^* \pi_1\big) 
\\   + 
\overline{\pr}_{[1,s + 2], *} 
\big(\overline{\pr}_{0}^*(K_C+2e) \cdot  \overline{\pr}_{k}^*e \cdot \prod_{l \neq j, k} \overline{\pr}_{0 l}^* \pi_1\big) \Big].
\end{multline*}
Consider the Cartesian square
\[
\begin{tikzcd}
    C\times C^{s+2}
    \arrow[r, "{\widehat{\overline{\pr}}_j}"]
    \arrow[d, "{\overline{\pr}_{[1,s+2]}}"']
    &
    C\times C^{s+1}
    \arrow[d, "{\overline{\pr}_{[1,s+1]}}"]
    \\
    C^{s+2}
    \arrow[r, "{\widehat{\overline{\pr}}_j}"]
    &
    C^{s+1}.
\end{tikzcd}
\]
By the {push-pull formula}, we have that
\begin{align*}
 \sum_{k \neq j} \overline{\pr}_{[1,s + 2], *} \big( \overline{\pr}_0^*(K_C + 4e) 
\cdot 
\prod_{l \neq j} \overline{\pr}_{0 l}^* \pi_1 \big) &  
 = (s+1) \widehat{\overline{\pr}}_j^*\, \overline{\pr}_{[1,s+1], *} \big(  \overline{\pr}_0^*(K_C + 4e) 
\cdot 
\prod_{l \neq j} \overline{\pr}_{0 l}^* \pi_1  \big) \\
 & =  (s+1)
 \widehat{\overline{\pr}}_{j}^* \beta_e^{s + 1}(K_C + 4e).
\end{align*}
where the last equality uses \Cref{betaz}. 
Similarly, we have 
\begin{equation*}
    \overline{\pr}_{[1,s + 2], *} \big( \overline{\pr}_0^*(K_C + 2e) \cdot \overline{\pr}_k^* e 
\underset{l \neq j, k}{\prod} \overline{\pr}_{0 l}^* \pi_1\big)= \overline{\pr}_k^* e \cdot \widehat{\overline{\pr}}_{j k}^* \beta_e^{s}(K_C + 2e).
\end{equation*}
Combining these three terms, we have the statement of the theorem.\qedhere
\end{proof}

\begin{corollary}\label{ceresatomoddiag}
    For $s\ge 0$, the projection of $(\phi \circ \sigma_{s + 2})^* \mathscr{F}([C]^e_{(s)})$ under $\pi_1^{\otimes (s + 2)}$ is $B^{s + 2}(C, e)$. 
    In particular if $[C]^e_{(s)}=0$, then $B^{s+2}(C,e)=\Gamma^{s+2}(C,e)=0$.
\end{corollary}
\begin{proof}
Let us consider equation \eqref{eq:equalityforFourier(C)}.
As $B^{s + 2}(C, e)$ is already in the 
subspace $\mathfrak{h}_{1\ldots1}(C^{s + 2})$, it is fixed by $\pi_1^{\otimes s+2}$. 
Note that all the other terms in \eqref{eq:equalityforFourier(C)} are of the from $C \times z'$ for some cycle $z'\in \CH_0(C^{s+1})$. 
Since $\pi_{1, *} (C) = 0$, by \eqref{decompositionofpushforward} all of the other terms vanish.
Therefore,
\[
\pi_{1,*}^{\otimes (s + 2)}\big((\phi \circ \sigma_{s + 2})^* \mathscr{F}([C]^e_{(s)}) \big) ={(-1)^s} B^{s + 2}(C, e).
\]
The rest of the statement follows from \Cref{thm:gammavanishiffB}.
\end{proof}

\begin{corollary}
We have that $\Gamma^n(C,e)=0$ for all $n> g+1$.
\end{corollary}

\subsection{From Diagonal Classes to Beauville Components}
We now show the converse of the results of the previous section, namely that the vanishing of diagonal classes imply the vanishing of some corresponding Beauville components. More precisely, we have the following theorem.

\begin{theorem}\label{moddiagtoceresa}
    Let $z\in \CH_i(C)$ and $s\ge 1$. Then we have
    \[
    {  \sigma_{s,*}(\gamma_e^s(z))_{(s-2i)} = s!\cdot \iota_{e,*}(z)_{(s-2i)} .}
    \]
    In particular, we have the implication 
    $$ \gamma_e^s(z)=0\quad \Longrightarrow \quad \iota_{e,*}(z)_{(s-2i)}=0. $$
\end{theorem}
\begin{proof}
We fix the value $n=s$ in \Cref{notationsection4}. Using the Cartesian diagram 
    \begin{center}
        \begin{tikzcd}
            C^s\times J \arrow[r,"\sigma_s\times \id"] \arrow[d,"\pr_{C^s}"] & J\times J \arrow[d,"\pr_{J, 1}"] 
            \arrow[r,"{\phi\times \id}"] & J^t\times J \arrow[d,"\pr_{J^t}"] \arrow[r,"\pr_{J}"] & J
            \\ C^s \arrow[r,"{\sigma_s}"] & J \arrow[r, "\phi"] & J^t,
        \end{tikzcd} 
    \end{center} 
by the {push-pull formula}, we have that
    \begin{align*}
    \mathscr{F}^{t}\big( (\phi\circ \sigma_s)_*\gamma_e^s(z) \big) &= \pr_{J,*}\big(\ch(\mathscr{P}^t)\cdot \pr^*_{J^t}(\phi\circ \sigma_s)_*(\gamma_e^s(z))\big)
    \\ &= \pr_{J,*}\big(\ch(\mathscr{P}^t)\cdot (\phi\times \id)_* (\sigma_s\times \id)_* \pr^*_{C^s}(\gamma_e^s(z))\big).
    \end{align*}
    By the projection formula with respect to $(\phi\times \id)\circ (\sigma_s\times \id)$ and the fact that $\pr_J\circ (\phi \times \id)\circ (\sigma_s\times \id))=\pr_J$, we get that
\begin{align*}
         \mathscr{F}^{t}\big((\phi\circ \sigma_s)_*\gamma_e^s(z) \big) &= \pr_{J,*}\big(\ch((\sigma_s\times \id)^*(\phi\times\id)^* \mathscr{P}^t)\cdot \pr^*_{C^s}(\gamma_e^s(z))\big)
         \\ &  = \pr_{J,*}\big(\ch((\sigma_s\times \id)^*\Lambda_\Theta)\cdot \pr^*_{C^s}(\gamma_e^s(z))\big),
\end{align*}
where the last equality uses the definition of $\Lambda_\Theta$. By \Cref{pullbackofmumford}, we have 
\[
(\sigma_s\times \id)^*\Lambda_{\Theta}=\sum_{j=1}^s\mathscr{L}_j,\ \textup{where} \ \mathscr{L}_j\colonequals( \overline{\pr}_{j}\times \id)^*(\mathscr{L})\ \textup{and} \ \mathscr{L}\colonequals (\iota_e\times \id)^*\Lambda_{\Theta}.
\]
Unwinding the definition of $\gamma_e^s(z)$, we get that
\begin{equation*}
\mathscr{F}^{t}\big((\phi\circ \sigma_s)_*\gamma_e^s(z) \big)
= 
\pr_{J,*}\big(
\prod_{j=1}^s \ch(\mathscr{L}_j)
\cdot \pr^*_{C^s}\pr_{[s+1,2s],*}\big( \pr_{[1,s]}^*\Delta_C^{(s)}(z) \cdot \prod_{j=1}^s \pr_{j\, j+s}^*(\pi_+) \big)\big).
\end{equation*}
By the Cartesian square
\[
\begin{tikzcd}
    C^{2s}\times J \arrow[rr,"{\pr_{[s+1,2s]}\times \id_J}"] \arrow[d,"{\pr_{C^{2s}}}"] & & C^s \times J \arrow[d,"{\pr_{C^s}}"]
    \\ C^{2s} \arrow[rr,"\pr_{[s+1,2s]}"] & & C^s,
\end{tikzcd}
\]
and the {push-pull formula}, we have 
\begin{equation*}
\mathscr{F}^{t}
\big((\phi\circ \sigma_s)_*\gamma_e^s(z) \big)
= 
\pr_{J,*}\Big(
\prod_{j=1}^s \ch(\mathscr{L}_j)
\cdot 
(\pr_{[s+1,2s]}\times \id_J)_* \pr^*_{C^{2s}}  \big( \pr_{[1,s]}^*\Delta_C^{(s)}(z) \cdot \prod_{j=1}^s \pr_{j\, j+s}^*(\pi_+) \big)\Big).
\end{equation*}
By the projection formula with respect to $\pr_{[s+1,2s]}\times \id_J$, we get
\begin{equation}\label{eq-Ftphi complicated}
\mathscr{F}^{t}\big((\phi\circ \sigma_s)_*\gamma_e^s(z) \big)= \pr_{J,*}
\Big(
\prod_{j=1}^s
\ch \big( (\pr_{s+j}\times \id)^*\mathscr{L} \big)
\cdot 
\pr_{[1,s]}^* \Delta^{(s)}_{C,*}(z)
\cdot
\prod_{j=1}^s \pr_{j\, j+s}^*(\pi_+) \Big).
\end{equation}
{We will compute $\mathscr{F}^t\big((\phi\circ \sigma_s)_*\gamma_e^s(z)_{(s-2i)} \big)$, which corresponds to the codimension $(s-i)$ part of \eqref{eq-Ftphi complicated}.}
We first note that if in the product $\prod_{j=1}^s
\ch ( (\pr_{s+j}\times \id_J)^*\mathscr{L} )$ there is a term with no contribution from {$
c_1 ((\pr_{s+k}\times \id_J)^*\mathscr{L})$ for some $1\leq k\leq s$}, then it can be written as $$\pr_{J,*}\left( (\widehat{\pr}_{k+s}\times \id_J)^*(z') \cdot \pr^*_{k\, k+s}(\pi_+) \right),$$ where $z'\in \CH(C^{2s-1}\times J)$. Therefore writing $\pr_J$ as $\pr_J\circ ( \widehat{\pr}_{k+s}\times \id_J)$, we can use projection formula with respect to $\widehat{\pr}_{k+s}\times \id_J$ and get
$$
\pr_{J,*}\left( 
z'\cdot 
\big(\widehat{\pr}_{k+s,*} \pr^*_{k\, k+s}(\pi_+)\big)\times J \right)
=
\pr_{J,*}\left( 
z'\cdot 
\big(\widetilde{\pr}_k^*\pr_{1,*}(\pi_+) \big)\times J \right),
$$
where $\widetilde{\pr}_k$ is the $k$-th projection $C^{2s-1}\rightarrow C$.
Since $\pr_{1,*}(\pi_+) =
[C]-\deg(e)\cdot[C] =0$, this term gives no contribution. {Furthermore any term with a higher power of $c_1 ((\pr_{s+k}\times \id_J)^*\mathscr{L})$ for some $1 \leq k \leq s$ will have codimension higher than $(s-i)$.} Therefore, 
$\mathscr{F}^t\left((\phi\circ \sigma_s)_*\gamma_e^s(z)_{({s-2i})} \right)$ 
is equal to
\begin{equation}\label{codimensionterm}
\pr_{J,*} \big(
\prod_{j=1}^s(\pr_{s+j}\times \id_J)^*c_1(\mathscr{L})\cdot \pr_{[1,s]}^*\Delta_C^{(s)}(z)  \prod_{j=1}^s \pr_{j\, j+s}^*(\pi_+)
\big). 
\end{equation}

All the terms in \eqref{codimensionterm} can be written as a pullback $(\widehat{\pr}_{1+s}\times \id_J)^{*}(z')$ with $z'\in \CH(C^{2s-1}\times~J)$ except for $(\pr_{s+1}\times \id_J)^*c_1(\mathscr{L})\cdot \pr^*_{1\, 1+s}(\pi_+)$: writing $\pr_J$ as $\pr_J\circ (\widehat{\pr}_{1+s}\times \id_J)$, we can use projection formula with respect to $\widehat{\pr}_{1+s}\times \id_J$ in \eqref{codimensionterm} and get
\begin{align*}
\pr_{J,*} \Big( z' \cdot (\widehat{\pr}_{s+1}\times \id_J)_*
& \big( (\pr_{s+1}\times \id_J)^*c_1(\mathscr{L})\cdot \pr^*_{1\, s+1}(\pi_+) \big) \Big) 
\\
= & \pr_{J,*} \Big( z' \cdot (\widehat{\pr}_{s+1}\times \id_J)_*\big( (\pr_{1\, s+1}\times \id_J)^*(C\times c_1(\mathscr{L}) \cdot \pi_+ \times J )  \big) \Big)
\\ 
= & \pr_{J,*}\Big( z' \cdot (\widetilde{\pr}_1\times \id_J)^*(\pr_1\times \id_J)_*(C\times c_1(\mathscr{L})\cdot \pi_+\times J) \Big),
\end{align*}
where the last equality comes from the {push-pull formula}.
Iterating, we see \eqref{codimensionterm} is equal to
$$ \pr_{J,*}\big(\Delta^{(s)}_{C,*}(z)\times J \cdot \prod_{j=1}^s (\overline{\pr}_j\times \id_J)^*(\pr_{C,1}\times \id_J)_*(C\times c_1(\mathscr{L})\cdot \pi_+\times J)\big),$$
where here $\pr_{C, 1}$ is the projection on the first component $C\times C\to C$. 
Writing $\pi_+=\Delta_C- (C\times e)$, we see that
\begin{align*}
& (\pr_{C, 1}\times \id)_*  (C\times c_1(\mathscr{L})\cdot \pi_+\times J) 
\\
& =  (\pr_{C, 1} \times \id)_*\big(C\times c_1(\mathscr{L})\cdot (\Delta_C\times \id_J)_*(C\times J) \big) - C\times \pr_{J,*}(c_1(\mathscr{L})\cdot e\times J) 
\\
& = (\Delta_C\times \id_J)^*(C\times c_1(\mathscr{L})) - C\times \pr_{J,*}(c_1(\mathscr{L})\cdot e\times J),
\end{align*}
where the last equality is given by projection formula with respect to $\Delta_C\times \id_J$.
The first term is equal to $c_1(\mathscr{L})$
and the term $\pr_{J,*}(c_1(\mathscr{L})\cdot e\times J)$ is the codimension one part of $\mathscr{F}^t((\phi\circ \iota_e)_*(e))$. By Fourier duality,  $\mathscr{F}^t$ induces an isomorphism $\CH_{0,(s)}(J^t)\xrightarrow{\sim}\CH^s_{(s)}(J)$. Therefore, we find that  $\pr_{J,*}(c_1(\mathcal{L})\cdot e\times J)= \mathscr{F}^t((\phi\circ \iota_e)_*(e)_{(1)})$. However, we have that $((\phi\circ \iota_e)_*(e))_{(1)}=\phi((\iota_{e,*}(e))_{(1)})=0$ since $(\iota_{e,*}(e))_{(1)}=\updelta^e_{(1)}=0$ by \Cref{yinI2}. Therefore the second term vanishes, and we get that \eqref{codimensionterm} is equal to
\begin{equation*}
\pr_{J,*}\big(\Delta_{C,*}^{(s)}(z)\times J \cdot \prod_{j=1}^s c_1(\mathscr{L}_j)\big).
\end{equation*}
Rewriting the diagonal term $\Delta_{C,*}^{(s)}(z)$ as $\overline{\pr}_1^*(z)\cdot \prod_{j=2}^{s}\overline{\pr}_{1 j}^*(\Delta_C)$, and noting that $c_1(\mathscr{L}_j)\cdot \overline{\pr}^*_{1 j}(\Delta_C)\times J= c_1(\mathscr{L}_1)\cdot \pr^*_{1 j}(\Delta_C)\times J$, we find that \eqref{codimensionterm} is equal to
\begin{align*}
    \pr_{J,*}\big(\overline{\pr}^*_1(z)\times J \cdot c_1(\mathscr{L}_1)^s\cdot \pr_{C^s}^*(\Delta^{(s)}_C)\big) &
    = \pr_{J,*}\big((\overline{\pr}_1\times \id)^*(z\times J\cdot c_1(\mathscr{L})^s)\cdot \pr_{C^s}^*(\Delta^{(s)}_C)\big)
    \\ &= \pr_{J,*}\big(z\times J \cdot c_1(\mathscr{L})^s\cdot ( {\overline{\pr}_1} \times \id)_*\pr_{C^s}^*(\Delta_C^{(s)})\big)
    \\ &= \pr_{J,*}(z\times J \cdot c_1(\mathscr{L})^s).
\end{align*}
This is precisely $s!$ times the codimension $(s-i)$ part of $\mathscr{F}^t((\phi\circ \iota_e)_*z)$, 
{which coincides with $\mathscr{F}^t(((\phi\circ \iota_e)_*z)_{(s-2i)})=\mathscr{F}^t(\phi_*( \iota_{e,*}(z)_{(s-2i)}))$. It therefore follows that }
\[
\mathscr{F}^{t}\big((\phi\circ \sigma_s)_*\gamma_e^s(z)_{(s-2i)} \big) = s!\cdot \mathscr{F}^t(\phi_*( \iota_{e,*}(z)_{(s-2i)})).
\] 
Applying $(\phi)^{-1}_* \circ (\mathscr{F}^t)^{-1}$ to both sides, we see $\sigma_{s,*}(\gamma_e^s(z))_{(s-2i)} = s! \cdot \iota_{e,*}(z)_{(s-2i)}$, as desired. 
\end{proof}

\begin{remark}
    The proof of this theorem actually gives a stronger statement than is necessary for our purposes. As the only terms from \eqref{eq-Ftphi complicated} that contribute are those where every embedding of $\mathscr{L}$ has a positive exponent, there is no contribution from terms of codimension $t < s - i$. In particular, $\sigma_{s, *}(\gamma_e^s(e))_{(t)} = 0$ for $t < s - 2i$.
\end{remark}

Besides establishing \Cref{maintheoremA}, we also get the following consequence  of  \Cref{ytogamma} and \Cref{moddiagtoceresa}.
\begin{corollary}\label{equivalentvanishingforzerocycles}
    If $s\ge 1$, then  $\updelta^e_{(s)}=0 \Longleftrightarrow \gamma^s_e(e)=0.$ 
\end{corollary}

Combining this result with the successive vanishing statements of \Cref{sec-successive vanishing} we also have the following.

\begin{corollary}\label{delta2vanish}
    If $\updelta^e_{(2)}=0$, then $\updelta^e=0$.
\end{corollary}
\begin{proof}
    Given the equivalence between the vanishing of $\updelta^e_{(s)}$ and that of $\gamma^s_e(e)$, the statement follows from \Cref{corollaryonvanishingofgamma}. 
    Alternatively, we can provide a different proof by observing that 
    \begin{align*}
    \sigma_{2, *}(\gamma^2_e(e)) & =  m_*\iota_{e,*}^2( \Delta_{C,*}(e)-e\times e)  = 
    [2]_* \iota_{e,*}(e) - \iota_{e,*}(e) \star \iota_{e,*}(e)
    \\
    & = [2]_*(\updelta^e+[0]) - (\updelta^e+[0])\star(\updelta^e+[0]) = [2]_*\updelta^e+ [0] - \updelta^e\star\updelta^e - 2\updelta^e -[0]
    \\
    & = [2]_* \updelta^e - 2\updelta^e + \updelta^e \star \updelta^e.
    \end{align*}
    By hypothesis, this cycle is zero.
    Now, in order to prove $\delta^e=0$, it is enough to show that $\updelta^e \in I^{\star n}$ for any $n\geq 1$ because $\bigcap_{n = 1}^{\infty} I^{\star n} = \{0\}$. 
    Let now $\updelta^e \in I^{\star n}$ for some $n \geq 2$. 
    Since $I^{\star n}= \oplus_{s \geq n} \CH_{0, (s)}(J)$, this means that $\delta^e_{(s)}=0$ for all $s<n$.
    Therefore 
    \[
    [2]_* \updelta^e - 2\updelta^e + \updelta^e \star \updelta^e \equiv 2^{n} \updelta^e - 2\updelta^e \equiv 0 \pmod{I^{\star (n + 1)}}
    \]
    which then implies $\updelta^e \in I^{\star (n + 1)}$.
\end{proof}

\begin{example}
    
    Let $E/\mathbb{C}$ be an elliptic curve and let $x,y \in E(\mathbb{C})$. Beauville and Voisin \cite[Lemma 2.5]{Beauville2001ONTC} showed that for $e=\frac{(x+y)}{2}$, we have
    $$
    4\cdot\gamma^2_{e}(e)=(x,x)-(x,y)-(y,x)+(y,y)=0 \quad \text{in }\CH^2(E^2;\mathbb{Z}).
    $$
    This can also be deduced from \Cref{equivalentvanishingforzerocycles} (and  Ro\u{\i}tman's Theorem), since $\updelta^e_{(s)}=0$ for any $s>g=1$. It is also shown in \cite[Proposition 4.2]{Beauville2001ONTC} that the same result holds for hyperelliptic curves $C$ when $y=w$ is a Weierstrass point. In our notation, this cycle is $\gamma^2_w(x)$, whose vanishing follows from that of $\Gamma^3(C,w)$ by \Cref{thm:vanishinglittlegamma}.
    We point out that 
    Qiu has generalized this fact for curves with real multiplication (see \cite[Section 2.6]{CQiu}).
\end{example}

\section{Applications}\label{sec:applications}
\subsection{New proof of Zhang's result}\label{sec:zhangnewproof} 
We briefly explain how our results recover \Cref{Zhangtheorem}, following an argument suggested to us by Ben Moonen.

By \Cref{ceresatomoddiag} and \Cref{moddiagtoceresa}, the vanishing of $[C]_{(1)}^{e}$ is equivalent to the vanishing of $\Gamma^3(C,e)$. Clearly, if $\Cer(C,e)=0$ then $[C]_{(1)}^e=0$ also vanishes. Thus it only remains to prove that $\Gamma^3(C,e)=0$ implies $\Cer(C,e)=0$, and that the vanishings of any of these cycles imply $e=\xi=K_C/(2g-2)$.

Assume that $\Gamma^3(C, e)$ vanishes, so that $B^3(C, e)$ also vanishes by \Cref{thm:gammavanishiffB}. Then by \Cref{vanishing of gamma n-12}, we have 
$$0=\gamma_e^1(K_C+2e) = K_C+2e-(2g)e,$$ 
which implies $e = \xi$. 
{Concretely, in the proof of this proposition, we view $B^3(C,e)$ as a correspondence between $C^2$ and $C^1$, and push $\Delta_C$ forward through this correspondence; a similar argument is used for instance in \cite[Proposition 3.1]{vial}}.  
Furthermore, the same proposition implies that $\gamma^2_e(e)=0$. Consequently, by \Cref{thm:recursiverelationbetagamma} we deduce that $\Gamma^n(C,e)=0$ for all $n\ge 3$. Finally in terms of the Beauville decomposition of the curve class, this means that $[C]_{(s)}^e=0$ for all $s\ge 1$, so that in particular $\Cer(C,e)=0$. Alternatively, we can recover the vanishing of the Ceresa cycle by computing $\sigma_{3,*}(\Gamma^3(C,e))$, using the fact that $\updelta^e=0$ by \Cref{equivalentvanishingforzerocycles} and \Cref{delta2vanish} (see \Cref{qiuremark}).

\begin{remark}\label{qiuremark}
After this article was first announced, Congling Qiu pointed out to us that our results fill a gap in Zhang's proof of \Cref{Zhangtheorem} (which carries over to a small error in the original statement of \cite[Theorem 1.5.5]{Zhang10}). Qiu had independently developed a similar argument in his note \cite{QiuNote}; we extend our gratitude to him for sharing a draft of that note. In summary, the statement of \cite[Theorem 1.5.5]{Zhang10} includes the following formula, which is stated unconditionally:
\begin{equation}\label{zhangformula}
    \sigma_{3,*}(\Gamma^3(C,e)) = \sum_{s\geq 0} (3^{2+s}-3\cdot 2^{2+s} +3)[C]^e_{(s)}.
\end{equation}
This formula is used in the proof of the theorem (see \cite[pg. 70]{Zhang10}): the equation \eqref{zhangformula} directly implies that if $\Gamma^3(C,e)=0$, then $[C]_{(s)}^e=0$ for all $s \geq 1$, as the coefficient $3^{2+s}-3\cdot 2^{2+s}+3$ is non-zero for $s \geq 1$. We note that by the definitions of $\sigma_3$ and $\Gamma^3(C,e)$ we have in fact
\begin{equation}\label{qiucorrect}
\sigma_{3,*}(\Gamma^3(C,e)) = [3]_* [C]^e - 3 \cdot [2]_* \big([C]^{e} \big)\star \iota_{e,*}(e) + 3 \cdot \big( [C]^e \big)\star \iota_{e,*}(e) \star \iota_{e,*}(e).
\end{equation}
(This equality is also used in our proof of \Cref{intModtoCer}). It is not clear whether \eqref{qiucorrect} necessarily implies \eqref{zhangformula} in general. However by \Cref{thm:vanishinglittlegamma}, \Cref{equivalentvanishingforzerocycles}, and \Cref{delta2vanish} there is a chain of implications
\[
\Gamma^3(C,e)=0 \ \  \Longrightarrow \ \  \gamma^2_e(e)=0 \ \ \Longrightarrow \ \  \updelta^e=0.
\]
The last vanishing is equivalent to $\iota_{e,*}(e)=[0]$; we remind the reader that $[0]\in \CH(J)$ is the identity element with respect to the Pontryagin product. Hence \eqref{qiucorrect} implies \eqref{zhangformula} at least under the assumption $\Gamma^3(C,e)=0$, or under the weaker assumption $\gamma_e^2(e)=0$ (which holds, for instance, when $e$ is rationally equivalent to a point). 
    
\end{remark}

\subsection{The $[C]_{(2)}$-component}\label{C2term}  
For this subsection, we work with the canonical choice of divisor $e=\xi=K_C/(2g-2)$. 

As we have noted, $\gamma^2_\xi(\xi)$ is a scalar multiple of the Faber--Pandharipande cycle. 
The work of Green and Griffiths \cite{GG03} (over $\mathbb{C}$) and Yin \cite{Yin15} (in arbitrary characteristic) shows that this cycle is generically non-vanishing for curves of genus $g \geq 4$. Combining this with our results, we may establish the nonvanishing of $[C]_{(2)}^{\xi}$ for $g\ge 4$ as an immediate corollary. We note that by the main theorem of \cite{CvG93}, $[C]_{(2)}^{\xi}$ vanishes algebraically for the generic curve of genus $4$.

\begin{proposition}\label{C2nonvanishing}
     For the generic curve $C$ of genus $g\ge 4$, we have $[C]^{\xi}_{(2)}\neq 0$. If $k$ is uncountable, then the same result holds for a very general curve $C$ with $g\ge 4$.
\end{proposition}
\begin{proof}
    By \Cref{maintheoremA}, the vanishing of $[C]^{\xi}_{(2)}$ is equivalent to that of $B^4(C,\xi)$, which in turn implies the vanishing of $\gamma_\xi^2(K_C + 2\xi) = 2g \cdot  \gamma_\xi^2(\xi)$ by the first part of \Cref{vanishing of gamma n-12}.  
    Hence the result follows from the main theorem of \cite{Yin15}.
\end{proof}

On the other hand, we have the opposite behaviour in genus 3. To prove this, we borrow some machinery from the work of Polishchuk \cite{Polish05}. For $e=x_0$ a $k$-point, Polishchuk studies the action of the Lie algebra $\mathfrak{sl}_2$ on $\CH(J)$ (see e.g. \cite[Section 9]{MoonenAWS}) and on its so-called \emph{tautological subring} $\mathcal{T}_{x_0}\CH(J)$. We recall that this subring contains the Beauville components $[C]^{x_0}_{(s)}$, and by \cite[Theorem 0.2]{Polish05} is generated by the classes $p_n\colonequals \mathscr{F}([C]^{x_0}_{(n - 1)})$ ($n \geq 1$) and $q_n \colonequals \mathscr{F}(\eta_{(n)})$ ($n \geq 0$), where $\eta\colonequals\iota_{x_0,*}(K_C/2)+[0]$. 
Although this work is restricted to the case when $e=x_0$ is a $k$-point, we can still recover sufficient information when $e=\xi$ is not necessarily a $k$-point.
\begin{theorem}\label{thm:c2vanishing}
Let $C$ be a genus 3 curve, then $[C]^\xi_{(2)} = 0$.
\end{theorem}
\begin{proof}
     We use the relation $[C]^{\xi} = [x_0 - \xi] \star [C]^{x_0}$. Taking the Fourier transform gives $\mathscr{F}([C]^{\xi}) = \mathscr{F}([x_0 - \xi]) \cdot \mathscr{F}([C]^{x_0})$. As remarked in \cite[\S 1]{Polish05}, we have $\mathscr{F}([x_0 - \xi]) = \exp(\Theta_{\xi-x_0} - \Theta)$ (where $\Theta_{\xi-x_0}$ and $\Theta$ are as in \Cref{subsectinoonthetadiv}). In particular, taking the codimension 3 component of these classes gives the following relation:
    $$\mathscr{F}([C]^{\xi}_{(2)}) = \mathscr{F}([C]^{x_0}_{(2)}) + (\Theta_{\xi-x_0} - \Theta) \mathscr{F}([C]^{x_0}_{(1)}) + \frac{1}{2} (\Theta_{\xi-x_0} - \Theta)^2 \mathscr{F}([C]^{x_0}_{(0)}).$$
    
    Since $\eta =\iota_{x_0,*}(K_C/2)+[0] $ maps to $2(\xi - x_0)$ under the Albanese map, the divisor $\Theta_{\xi-x_0} - \Theta$ is equal to $-\frac{1}{2}\mathscr{F}(\eta_{(1)})$ (see \cite[\S 1]{Polish05}). We thus obtain the following expression for $\mathscr{F}([C]^{\xi}_{(2)})$ in terms of $p_n$ and $q_n$:
    $$8\mathscr{F}([C]^{\xi}_{(2)}) = 8p_3 - 4q_1 p_2 + q_1^2 p_1.$$
    This expression can be simplified using the relations in the tautological ring to give the vanishing of the right hand side. Explicitly, as in \cite[Proposition 2.1, Corollary 2.2 (ii)]{Yin15}, we have $4q_2 - q_1^2 = 0$. Moreover, by \cite[Theorem 0.2]{Polish05}, we have $\mathbf{f}(p_2^2 - 2p_1 q_2 p_2 + 2p_1 p_3) = 8p_3 - 4p_2 q_1 - 4q_2 p_1 + 2p_1 q_1^2$, where $\mathbf{f}$ denotes the usual generator for the Lie algebra $\mathfrak{sl}_2$. But the cycle $p_2^2 - 2p_1 q_2 p_2 + 2p_1 p_3$ vanishes for dimension reasons, so $[C]^{\xi}_{(2)}=0$.
\end{proof}

\begin{remark}
Under the assumption $[C]^\xi_{(2)}=0$, and working with the cycle $e=\xi$, we can strengthen some of the successive vanishing results of \Cref{moddiagcycles}.
Indeed the vanishing of $B^4(C, \xi)$ implies the vanishing of $\gamma^2_\xi(\xi)$. 
But as we saw in \Cref{successivegammavanishing} and \Cref{equivalentvanishingforzerocycles}, the vanishing of $\gamma^2_\xi(\xi)$ implies the vanishing of $\gamma^n_\xi(\xi)$ for all $n \geq 2$, and that of $\updelta^\xi$. By \Cref{thm:recursiverelationbetagamma}, this also implies $B^5(C, \xi) = 0$, and so $[C]^\xi = [C]^{\xi}_{(0)} + [C]^{\xi}_{(1)}$. 

For general $e$, the situation is slightly murkier. 
By \Cref{thm:vanishinglittlegamma}, we have $\gamma^n_e(e) = 0$ for all $n \geq 3$, and so $\updelta^e$ has at most one non-trivial Beauville component ($s = 2$). We also get $\Gamma^n(C, e) = 0$ for all $n \geq 6$ by \Cref{vanishinglargegamma}, which forces $[C]^e = [C]^e_{(0)} + [C]^e_{(1)} + [C]^e_{(3)}$.  
\end{remark}

\begin{remark}
In forthcoming work, we will generalize Polishchuk's results on the tautological subring $\mathcal{T}_{x_0} \CH(J)$ to the case where the basepoint $x_0$ is replaced by an arbitrary divisor $e$.
In particular, we compute an explicit set of generators and relations for $\mathcal{T}_e \CH(J)$. 
The embedding with respect to the basepoint $x_0$ was sufficient for the proof of \Cref{thm:c2vanishing}, since the computations simplify considerably for curves of small genus.
\end{remark}

\subsection{Integral aspects}\label{sec:integral}
In this subsection, we prove the integral refinement \Cref{integrality} for Zhang's result. We combine our methods with results of Moonen and Polishchuk \cite{MP10}, who introduce a grading for the ring $\CH(J;\bbZ)$ satisfying properties similar to those of the usual Beauville decomposition. 

We let $\xi_{\mathrm{int}}\in \CH_0(C;\bbZ)$ be an integral representative for $\xi$, so that there exists an integer $N\neq 0$ satisfying $N(2g-2)\cdot\xi_{\mathrm{int}} = N\cdot K_C$. 
We fix an integer $d\in \mathbb{Z}$.

\begin{proposition}\label{torsionorderofGamma3fromCeresa}
If $\Cer(C,\xi_{\mathrm{int}})\in \CH_1(J;\bbZ)[d]$, then $\Gamma^3(C,\xi_{\mathrm{int}})\in \CH_1(C^3;\bbZ)[2\cdot d]$.
\end{proposition}
\begin{proof}
Motivated by the proof of the corresponding statement with rational coefficients given in \Cref{pullback of fourier transform} and \Cref{expression of fourier transform}, we are led to compute the cycle 
\begin{equation}\label{eqtocompute}
    \pr_{C^3,*}\big(\pr^*_J\Cer(C,\xi_{\mathrm{int}})\cdot\prod_{j=1}^3 c_1(\mathscr{L}_j) \big),
\end{equation}
where $\mathscr{L}_j=(\id_J\times (\iota_{\xi_{\mathrm{int}}}\circ \pr_j))^*\Lambda_{\Theta}$ is a line bundle on $J\times C^3$. Given that $([-1]\times \id_J)^*\Lambda_{\Theta}=-\Lambda_{\Theta}$, it follows that \eqref{eqtocompute} is equal to 
\begin{equation*}
    2 \pr_{C^3,*}\big( \pr^*_J [C]^{\xi_{\mathrm{int}}}\cdot\prod_{j=1}^3 c_1(\mathscr{L}_j) \big).
\end{equation*}
By the Cartesian diagram from the proof of \Cref{expression of fourier transform}, we find that 
\begin{align*}
    2 \pr_{C^3,*}\big( \pr^*_J [C]^{\xi_{\mathrm{int}}} \cdot \prod_{j=1}^3 c_1(\mathscr{L}_j)\big)= \  &2 \pr_{C^3,*}\Big((\iota_{\xi_{\mathrm{int}}}\times \id_{C^3})_*[C\times C^3]\cdot\prod_{j=1}^3 c_1(\mathscr{L}_j) \Big).\\ 
    \intertext{By the projection formula this simplifies to}
    &   2\pr_{C^3,*}\big(\prod_{j=1}^3(\iota_{\xi_{\mathrm{int}}}\times \id_{C^3})^*c_1(\mathscr{L}_j)\big),\\
    \intertext{and by \Cref{lem-pullbackMB} this equals}
    & 2\pr_{C^3,*}\big(\prod_{j=1}^3 \overline{\pr}_{0j}^*\pi_1\big).
\end{align*}
Finally by \Cref{betaz}, this is equal to $2 B^3(C,\xi_{\mathrm{int}})$. 
The equivalence between the vanishing of $B^3(C,\xi_\mathrm{int})$ and $\Gamma^3(C,\xi_\mathrm{int})$ remains, since $B^3(C, \xi_\mathrm{int})$ is torsion, and by the proof of \Cref{vanishing of gamma n-12}, $\gamma^2_{\xi_\mathrm{int}}(\xi_\mathrm{int})$ is also torsion. However, as $\gamma^2_{\xi_\mathrm{int}}(\xi_\mathrm{int})$ is in the kernel of the Albanese map, by Ro\u{\i}tman's Theorem \cite{Rojtman} (completed by Milne in positive characteristic \cite{MilneZC}), this cycle vanishes integrally.
\end{proof}

\begin{remark}
We note that \Cref{torsionorderofGamma3fromCeresa} is the optimal result one can expect at this level of generality. 
Indeed, Moonen and Polishchuk provided an example of a hyperelliptic curve $C$ (so that $\Cer(C,\xi_{\mathrm{int}})=0$) for which $\Gamma^3(C,\xi_{\mathrm{int}})$ is non-trivial and 2-torsion (see \cite[Remark 2.6]{MP_DivPow}).
\end{remark}

\begin{remark}
In the case of hyperelliptic curves, \Cref{torsionorderofGamma3fromCeresa} recovers a result of Moonen and Polishchuk \cite[Proposition 2.5]{MP_DivPow} stating that $\Gamma^3(C,x_0)$ is $2$-torsion for a Weierstrass point $x_0\in C(k)$. This is in turn a strengthening of a result of Gross and Schoen \cite[Corollary 4.9]{GS95}, who show that $\Gamma^3(C, x_0)$ is 6-torsion modulo algebraic equivalence.
\end{remark}

\begin{remark}\label{GammaandB3havesameorder}
    Let us highlight that the last argument in the proof of \Cref{torsionorderofGamma3fromCeresa} shows that the order of $B^3(C,e)$ coincides with the one of $\Gamma^3(C,e)$ in $\CH_1(C^3; \mathbb{Z})$. 
    {Moreover, if $\Gamma^3(C,\xi_{\mathrm{int}})$ is $d$-torsion then $\gamma^2_{{\xi_\mathrm{int}}}(\xi_{\mathrm{int}})=0$ by the proof of \Cref{torsionorderofGamma3fromCeresa}. 
    Applying \Cref{thm:recursiverelationbetagamma} (which holds integrally) with $m=1$ and $z=d\cdot [C]$, and arguing by induction, yields that $\Gamma^n(C,\xi_{\mathrm{int}})$ is $d$-torsion for $n\ge 3$. }
\end{remark}

\begin{theorem}\label{intModtoCer}
Let {$M_{g+1} = \prod_{\text{prime }p \leq g+1} p^{\ell_p}
$}, with $\ell_p=\lfloor \frac{g}{p-1}\rfloor$ if $p\ge 3$, and $\ell_2=g-1$.
If $\Gamma^3(C,\xi_{\mathrm{int}})\in \CH_1(C^3;\bbZ)[d]$, then $\Cer(C,\xi_{\mathrm{int}})\in \CH_1(J;\bbZ)[M_{g+1}\cdot d]$.
\end{theorem}
\begin{proof}
By \cite[Theorem 4]{MP10}, the group $\mathrm{CH}_{i}( J;\bbZ)$ admits a direct sum decomposition
\begin{equation}\label{MPdecomp}
\CH_i(J;\bbZ)=\bigoplus_{m=0}^{g+i}\CH^{[m]}_i(J;\bbZ)
\end{equation}
with the following property: for all $z\in \CH^{[m]}_i(J;\bbZ)$ and all $n\in \mathbb{Z}$, 
\begin{equation}\label{eq-mpproperty}
    [n]_*(z)=n^m\cdot z +y,\quad y\in \bigoplus_{m'>m}\CH^{[m']}_i(J;\bbZ).
\end{equation}  
It follows that the pushforward $[n]_{*}$ respects the filtration $\mathrm{Fil}^{m}\CH_{i}( J ;\bbZ) = \CH_i^{[\geq m]}( J;\bbZ)$, and furthermore, $[n]_{*}$ acts on $\mathrm{gr}_{\mathrm{Fil}}^{m}\mathrm{CH}_{i}( J ;\bbZ)$ by multiplication by $n ^{m}$. This is analogous to the action of $[n]_{*}$ on the Beauville decomposition (which is only defined over $\mathbb{Q}$-coefficients), though the decompositions are in fact distinct even with $\Q $-coefficients. We note also that the indexing is shifted by $2i$.

By \cite[\S7]{MP10} and by \cite[Corollary 3.8(i)]{MP10}, for dimension $i=1$ we can start the decomposition \eqref{MPdecomp} at $m=2$ and write
$$ [C]^{\xi_{\mathrm{int}}}=[C]^{[2]}+\cdots+[C]^{[g+1]}. $$
It then follows from \eqref{eq-mpproperty} that when the decomposition \eqref{MPdecomp} is applied to the cycle $ \Cer(C, \xi_{\mathrm{int}})$, it starts at $m=3$:
$$ \Cer(C,\xi_{\mathrm{int}})=(\mathrm{id}- [-1])_*[C]^{\xi_{\mathrm{int}}} =  \Cer(C,\xi_{\mathrm{int}})^{[3]}+\cdots+ \Cer(C,\xi_{\mathrm{int}})^{[g+1]} .$$
We stress that the equality $\Cer(C,\xi_{\mathrm{int}})^{[m]}=(1-(-1)^m)\cdot[C]^{[m]}$ does not necessarily hold.

Assume that $\Gamma^3(C,\xi_{\mathrm{int}})$ is $d$-torsion.
By \Cref{GammaandB3havesameorder} $\Gamma^n(C,\xi_{\mathrm{int}})$ is $d$-torsion for $n\ge 3$.
Since $\updelta^{\xi_{\mathrm{int}}}=\iota_{\xi_{\mathrm{int}},*}(\xi_{\mathrm{int}})-[0]$ is torsion and in the kernel of the Albanese, it vanishes integrally. 
Hence $\iota_{\xi_{\mathrm{int}},*}(\xi_{\mathrm{int}})=[0]$. {Applying $(2)$ of \Cref{explicitexpression}, we find that}
\begin{align}
    0= d\cdot \sigma_{n,*}(\Gamma^n(C,\xi_{\mathrm{int}})) 
    &= d\cdot \sum_{k=0}^{n-1} (-1)^{k}\binom{n}{k}[n-k]_*([C]^{\xi_\mathrm{int}})\star \iota_{\xi_{\mathrm{int}},*}(\xi_{\mathrm{int}})^{\star k}\nonumber\\
    & = d\cdot \sum_{k=0}^{n-1} (-1)^{k}\binom{n}{k}[n-k]_*([C]^{\xi_\mathrm{int}})\label{eq-integralbinomialvanishing},
\end{align}
where the second equality uses that $\iota_{\xi_{\mathrm{int}},*}(\xi_{\mathrm{int}}) = [0]$ is the identity with respect to the Pontryagin product. We take the pushforward of \eqref{eq-integralbinomialvanishing} by $(\mathrm{id}-[-1])$ and get a similar result for $\Cer(C, \xi_{\mathrm{int}})$:
\begin{equation}\label{eq-intceresaadded}
    0 = d\cdot \sum_{k=0}^{n-1} (-1)^{k}\binom{n}{k}[n-k]_*(\Cer(C, \xi_{\mathrm{int}})).
\end{equation}
We conclude that the $\CH^{[3]}_1(J)$ component of $(\id-[-1]_*)\big(d\cdot\sigma_{3,*}(\Gamma^3(C,\xi_{\mathrm{int}}))\big)$ is given by
\begin{equation}\label{eq-base induction int}
 \big(d\cdot\sum_{k=0}^{2}(-1)^k\binom{3}{k}(3-k)^3 \big)\cdot\Cer(C,\xi_{\mathrm{int}})^{[3]}=6d \cdot\Cer(C,\xi_{\mathrm{int}})^{[3]}=0.
\end{equation}

For each integer $s \geq 3$, denote by $\mathcal{P}_{s}$ the following set of primes 
\begin{equation*}
\mathcal{P}_{s} = \{ p \leq s : (p-1) \mid (s-1) \}.
\end{equation*}
For $m \geq 3$, let $M_{m}$ be the positive integer $\prod_{s=3}^{m}\prod_{p \in \mathcal{P}_{s}}^{}p$. We note that $v_{p}\left( M_{m} \right) =\lfloor \frac{m-1}{p-1}\rfloor$ for primes $p\ge 3$ and $v_2\left( M_{m} \right)  = m-2$. Of course also $M_{m'}\mid M_{m}$ for $m' \leq m$. 

We claim that for all $m \geq 3$, we have $dM_m \cdot \mathrm{Cer}( C, \xi _{\mathrm{int}} ) ^{[m]}=0$. We proceed by induction on $m$; the equation \eqref{eq-base induction int} implies the claim is true for $m=3$. Assume the claim holds for all $m' \leq m$.
Applying the decomposition \eqref{MPdecomp} to the cycle $dM_m \cdot \mathrm{Cer}( C, \xi _{\mathrm{int}} )$, and using $M_{m'}\mid M_m$ for $m'\le m$, we find that
$$ dM_m \cdot \mathrm{Cer}( C, \xi _{\mathrm{int}} )= (dM_m \cdot \mathrm{Cer}( C, \xi _{\mathrm{int}} ))^{[m+1]}+(dM_m \cdot \mathrm{Cer}( C, \xi _{\mathrm{int}} ))^{[m+2]}+\cdots, $$ 
with $(dM_m \cdot \mathrm{Cer}( C, \xi _{\mathrm{int}} ))^{[m+1]}=dM_m \cdot (\mathrm{Cer}( C, \xi _{\mathrm{int}} )^{[m+1]})$ by \eqref{eq-mpproperty}.

Suppose that $m+1$ is even, then we have $\mathcal{P}_{m+1}=\{ 2 \}$ and $M_{m+1} = 2M_m$. We have \[
	\big([-1]_{*}  \Cer\left( C,\xi _{\mathrm{int}} \right)\big)^{[m+1]}  = \Cer\left( C,\xi _{\mathrm{int}} \right)^{[m+1]} + \big(  [-1]_{*}\Cer( C,\xi _{\mathrm{int}} )^{[\leq m]}\big) ^{[m+1]}
.\]
Since $[-1]_{*}  \Cer( C,\xi _{\mathrm{int}})=-\Cer( C,\xi _{\mathrm{int}})$, we find that $dM_{m}\cdot \Cer\left( C,\xi _{\mathrm{int}} \right) ^{[m+1]}$ is 2-torsion, and so $dM_{m+1}\cdot \Cer\left( C,\xi _{\mathrm{int}} \right)^{[m+1]}=~0$. 

Now suppose that $m+1$ is odd. Taking the 
$\CH^{[m+1]}_1(J)$-component of equation \eqref{eq-intceresaadded} multiplied by $M_m$ yields
\begin{equation*}
\big(dM_{m}\cdot\sum_{k=0}^{n-1}(-1)^k\binom{n}{k}(n-k)^{m+1} \big)\cdot\Cer(C,\xi_{\mathrm{int}})^{[m+1]}=0.
\end{equation*}
Note that the coefficient is non-zero if and only if $n\le m+1$ (see \cite[Lemma 4.10]{Vois15}). Using \Cref{gcd}, we deduce that $dM_{m+1}\cdot\Cer(C,\xi_{\mathrm{int}})^{[m+1]}=0$.
\end{proof}

\begin{remark}
    The constants $M_m$ are closely related to the Hirzebruch numbers $T_m$, which appear as the denominators of the degree $m$ component of the Todd power series (see \cite{Hirz}, Lemma 1.7.3). More precisely, we have $2M_m=T_{m-1}$. These constants play a role in the integral formulation of the Grothendieck--Riemann--Roch Theorem (see \cite{Pappas}), and also feature in the recent construction of an integral Fourier transform \cite{integralFourier}. This is suggestive that this might be the limit of our methods.
\end{remark}

\begin{lemma}\label{gcd}
Given $s\ge 3$ odd, the greatest common divisor of all \[f(n)= \sum_{k=0}^{n-1} (-1)^{k}\binom{n}{k} (n-~k)^s\] for $3\le n \le s$ is given by the product $\prod_{p\in \mathcal{P}_s}p$, where $ \mathcal{P}_s=\{p \le s \text{ prime } : \  (p-1)\mid (s-1)\}. $
\end{lemma}
\begin{proof}
{We will use the identity
   \begin{equation}\label{eq-differenceoperator}
\sum_{k=0}^{n-1}(-1)^k\binom{n}{k}(n-k)^{t} =0, \quad 1 \leq t \leq n-1.
   \end{equation}
Indeed this expression is $D^{n}(x^{t})\mid_{x=0}$, for $D$ the difference operator $D(h)(x) = h(x+1) - h(x)$.}

First, assume $p$ is a prime with $(p - 1) \mid (s - 1)$. Then, $(n - k)^s \equiv (n - k) \pmod{p}$. 
This implies in paticular that $f(n) \equiv \sum_{k = 0}^{n - 1} (-1)^k {n \choose k} (n - k) = 0 \pmod{p}$, where we use \eqref{eq-differenceoperator} with $t=1$. Hence every $p \in \mathcal{P}_s$ divides $f(n)$ for each  $3 \leq n \leq s$. 

   Now assume $p$ is an odd prime dividing $f(n)$ for all $3 \leq n \leq s$. In particular, $p$ divides \[f(4) + 4f(3) = 4^s - 6 \cdot 2^s - 8 = (2^s - 2)(2^s - 4).\] 
   Therefore, the integer $\mathrm{ord}_{p}\left( 2 \right) $ (the multiplicative order of $2$ modulo $p$) divides either $s - 1$ or $s - 2$. {In the first case, we show inductively that $\mathrm{ord}_{p}\left( i \right) \mid (s-1)$ for all $3 \leq i \leq p - 1$; in the second case we show $\mathrm{ord}_p(i) \mid (s-2)$ for all $3 \leq i \leq p-1$. Indeed in the first case, suppose $\mathrm{ord}_p(j) \mid (s-1)$ holds for all $3 \leq j \leq i-1$. Since $p$ divides $f(i)$, we have 
   \begin{equation}\label{eq-inductive assumption}
   i^s + \sum_{k = 1}^{i - 1} (-1)^k {i \choose k}(i - k)^s \equiv 0 \pmod{p} .
   \end{equation}
 Taking \eqref{eq-differenceoperator} in the case $t=1$, we see
\[ i = -\sum_{k = 1}^{i - 1} (-1)^k {i \choose k} (i - k) \equiv i ^{s} \pmod{p},\] 
and so $i^{s-1}\equiv 1 \pmod{p}$, as desired. The second case follows analogously: if $\mathrm{ord}_p(j) \mid (s-2)$ holds for all $3 \leq j \leq i-1$, then taking \eqref{eq-differenceoperator} for $t=2$, we have
\[ i^2 = -\sum_{k = 1}^{i - 1} (-1)^k {i \choose k} (i - k)^2 \equiv i ^{s} \pmod{p}.\] 
If $i^{s - 1} \equiv 1$ (resp. $i^{s - 2} \equiv 1$) holds modulo $p$ for all $1 \leq i \leq p -1$, it follows that $(p - 1) \mid (s - 1)$ (resp. $\left( p-1 \right) \mid \left( s-2 \right) $). As $s$ and $p$ are odd, we must be in the first case, and so $p \in \mathcal{P}_s$.}

Finally, we show that the greatest common divisor is square-free. At the prime $2$ we simply note that $f(3) \equiv 3^s + 3 \equiv 2 \pmod{4}$. For odd primes $p \in \mathcal{P}_{s}$, the above argument showing $i^{s - 1} \equiv 1$, also holds modulo $p^2$. Then, we see \[f(p) = p^s + \sum_{k = 1}^{p - 1} (-1)^k {p \choose k} (n - k)^s \equiv \sum_{k = 1}^{p - 1} (-1)^k {p \choose k} (n - k) \equiv  -p \pmod{p ^2}.\]  
Hence the greatest common divisor is not divisible by $p ^2$ for any $p \in \mathcal{P}_{s}$, as desired.
\end{proof}

Let us recall the vanishing criterion of \cite{QZ24}: if a curve $C$ is endowed with the action of a group $G$ by automorphisms such that $(H^1(C)^{\otimes 3})^G=0$, then $\Gamma^3(C,\xi)$ vanishes in $\CH^2(C^3)$. 
This criterion admits an \emph{integral refinement}, thereby providing an upper bound on the order of the modified diagonal whenever it applies.

\begin{theorem}\label{thm:integralversofQZ24}
Let $G$ be a subgroup of $Aut(C)$ such that $(H^1(C)^{\otimes 3})^G =0$, where $H^1(C)$ is any Weil cohomology theory on which $\End(J)$ acts faithfully (for instance, Betti cohomology in characteristic 0). Then we have the following torsion bound,
\[
N\cdot (2g-2)|G|\cdot \Gamma^3(C,\xi_{\mathrm{int}})=0 \textup{ in } \CH^2(C^3; \mathbb{Z}),
\]
where as before $N$ is such that $N(2g - 2)\cdot\xi_{\mathrm{int}} = N\cdot K_C$.
\end{theorem}

\begin{proof}
    Let us check that all the steps in the proof of \cite[Theorem 1.2.1]{QZ24} remain valid integrally:
    \begin{enumerate}
        \item[1.] Viewing elements of $\CH_1(C^2;\bbZ)$ as self correspondences on $C$, we have an exact sequence $$0\to \CH_0(C;\bbZ)\oplus \CH_0(C;\bbZ)\xrightarrow{\pr_1^*\oplus \pr_2^*} \CH_1(C^2;\bbZ)\to \End(J)\to 0.$$
        The section $f\mapsto (\iota_{\xi_{\mathrm{int}}}^2)^*c_1((f\times \id_J)^*\Lambda_{\Theta})$ induces an isomorphism $\End(J)\xrightarrow{\sim}\pi_{1,*}\CH_1(C^2;\bbZ)$, where we recall that $\pi_1=\Delta_C-(C\times \xi_{\mathrm{int}})-(\xi_{\mathrm{int}}\times C)$.
        \item[2.] For $z\in [\pi_{1,*}\CH_1(C^2;\bbZ)]^{\otimes n}$, let $\mathrm{cl}(z)\in \End(H^1(C)^{\otimes n})$ be the corresponding endomorphism by functoriality. Since $\End(J)$ acts faithfully on $H^1(C)$, we have that $\mathrm{cl}(z)=0$ if and only if $z=0$.
        \item[3.] For $g\in G$, let $\Gamma_g\in \CH_1(C^2;\bbZ)$ be the graph of the automorphism of $C$ defined by $g$. Given that $N\cdot K_C=N(2g-2)\cdot\xi_{\mathrm{int}}$ is $G$-invariant, we have that 
        $$
        N(2g-2)\cdot \Gamma_{g,*} \circ \pi_{1,*}=N(2g-2)\cdot \pi_{1,*}\circ \Gamma_{g,*}.
        $$ 
        Consider the correspondence $z=\sum_{g\in G}(\Gamma_{g,*}\circ \pi_{1,*})^{\otimes 3}$. Since $\Gamma_{g,*}^{\otimes 3}(\Delta_C^{(3)})=\Delta_C^{(3)}$, we deduce that $N\cdot (2g-2)\cdot z(\Delta_C^{(3)})=N\cdot(2g-2)|G|\cdot B^3(C,\xi_{\mathrm{int}})$. 
        On the other hand, by definition of $z$ and the hypothesis $(H^1(C)^{\otimes 3})^{G}=0$, we find that $\mathrm{cl}(z)$ acts trivially on $H^1(C)^{\otimes 3}$, hence $z=0$. 
        Therefore, we get that $N\cdot(2g-2)|G|\cdot B^3(C,\xi_{\mathrm{int}})=0$. 
    \end{enumerate}
    One concludes by \Cref{GammaandB3havesameorder}.
\end{proof}

\begin{example}
Consider the bielliptic Picard curve $C: y^3=x^4+1$ which has automorphism group of order $48$, and for which the criterion applies (\cite[Example 3.1]{LS24}). 
Then by \Cref{thm:integralversofQZ24}, the order of $\Gamma^3(C,\infty)$ divides $ (2\cdot 3 -2) \cdot 48 =  192$, where $\infty$ is the unique point at infinity. 
So by \Cref{intModtoCer}, $\Cer(C,\infty)$ has order dividing $2304$. We remark that these bounds hold for any divisor $e$ with $4e = K_C$. One may expect $B^3(C, \infty)$ to have a smaller order than $B^3(C, e)$ for a general $e$ of this type, but we cannot prove this at present.
\end{example}

\begin{example}\label{hypergeom_nonvan}
    Let $n$ be a positive integer. We consider the hypergeometric curves $C_{n, \lambda} : (x^n - 1)(y^n - 1) = \lambda$, for a parameter $\lambda \in \mathbb{P}^1_k \setminus \{0, 1, \infty\}$, as defined by M. Asakura and N. Otsubo \cite[Section 2.1]{AsaOtsu}. In \cite{esknem}, P Eskandari and Y Nemoto show that for $p \geq 5$ a prime, and $e$ a ``cusp'' (one of the points at infinity), $\Gamma^3(C_{p, \lambda}, e)$ is not $l$-torsion for any positive integer $l \leq \frac{p - 1}{2}$. Combining this with \Cref{torsionorderofGamma3fromCeresa} shows that $\Cer(C_{p, \lambda}, e)$ is not $l$-torsion for any $l \leq \frac{p - 1}{4}$. In particular as $l=1$ is in this range, the Ceresa cycle is non-zero for this choice of basepoint. 
    
    On the other hand when $p = 3$, Eskandari and Nemoto show that $\Gamma^3(C_{3,\lambda},e)$ is torsion for every $\lambda$ and every choice of basepoint $e$. We note that the curves $(x^3 - 1)(y^3 - 1) = \lambda$ are a 1-dimensional family of genus 4 curves with automorphism group containing $S_3^2$; indeed the automorphism group contains $(x, y) \mapsto (\omega x, \omega y)$, $(x, y) \mapsto (\frac{\mu}{x}, \frac{\mu}{y})$, $(x, y) \mapsto (\omega x, \omega^{-1} y)$ and $(x, y) \mapsto (y, x)$, where $\omega$ is a cube root of unity, and $\mu^3 = 1 - \lambda$. This explicit family therefore realises the genus 4 family from \cite{QZ24}. Then by \Cref{thm:integralversofQZ24}, if $e$ is any divisor satisfying $(2\cdot 4 - 2)e = K_{C_{3,\lambda}}$, the cycle $\Gamma^3(C_{3,\lambda}, e)$ has order dividing $(2 \cdot 4 - 2) \cdot 36 = 216$, and so by \Cref{intModtoCer} the cycle $\Cer(C_{3,\lambda}, e)$ has order dividing $216\cdot360 = 77760$.
\end{example}

\bibliographystyle{alpha}
\bibliography{literature}
\end{document}